\documentclass[onefignum,onetabnum]{siamart171218}

\usepackage{graphicx}
\usepackage{booktabs}
\usepackage{amsmath}
\usepackage[english]{babel}
\usepackage[utf8]{inputenc}
\usepackage{amssymb}
\usepackage[nocompress]{cite}
\usepackage{hyperref}
\usepackage{multirow}
\usepackage{bbm}
\usepackage{algorithm}
\usepackage{enumitem}
\usepackage{algpseudocode}

\usepackage{tikz}
\usetikzlibrary{arrows.meta, positioning, calc}

\makeatletter
\newcommand{\resulttablecaption}{\def\@captype{table}\caption}
\AtBeginDocument{
  \def\refstepcounter@optarg[#1]#2{%
    \cref@old@refstepcounter{#2}%
    \cref@constructprefix{#2}{\cref@result}%
    \@ifundefined{cref@#1@alias}%
      {\def\@tempa{#1}}%
      {\def\@tempa{\csname cref@#1@alias\endcsname}}%
    \protected@edef\cref@currentlabel{%
      [\@tempa][\arabic{#2}][\cref@result]%
      \csname p@#2\endcsname\csname the#2\endcsname}}%
}
\makeatother

\newcommand{\bA}{{\bf A}}
\newcommand{\bB}{{\bf B}}
\newcommand{\bC}{{\bf C}}
\newcommand{\bD}{{\bf D}}

\newcommand{\bF}{{\bf F}}
\newcommand{\bG}{{\bf G}}
\newcommand{\bH}{{\bf H}}
\newcommand{\bI}{{\bf I}}
\newcommand{\bJ}{{\bf J}}
\newcommand{\bL}{{\bf L}}
\newcommand{\bM}{{\bf M}}

\newcommand{\bP}{{\bf P}}
\newcommand{\bQ}{{\bf Q}}

\newcommand{\bS}{{\bf S}}

\newcommand{\ba}{{\bf a}}
\newcommand{\bb}{{\bf b}}
\newcommand{\bc}{{\bf c}}

\newcommand{\bff}{{\bf f}}
\newcommand{\bg}{{\bf g}}

\newcommand{\bp}{{\bf p}}

\newcommand{\bx}{{\bf x}}
\newcommand{\by}{{\bf y}}

\newcommand{\bzero}{{\bf0}}

\newcommand{\bmu}{{\bf \mu}}

\newcommand{\bepsilon}{{\boldsymbol{\epsilon}}}

\newcommand{\calN}{\mathcal{N}}

\newcommand{\calR}{\mathcal{R}}

\newcommand{\calY}{\mathcal{Y}}

\newcommand{\bbR}{\mathbb{R}}

\newsiamremark{remark}{Remark}
\makeatletter
\providecommand*{\toclevel@remark}{3}
\providecommand*{\toclevel@proposition}{3}
\providecommand*{\toclevel@algorithm}{1}
\makeatother

\DeclareMathOperator*{\argmin}{arg\,min}

\newcommand{\TheTitle}{Marginalization-Based Regularization for Separable Nonlinear Least Squares Problems with Applications to  Semi-Blind  Deblurring}
\newcommand{\ShortTitle}{Marginalization-Based Regularization}
\newcommand{\TheAuthors}{J. T. Dworaczyk and M. I. Espa\~nol}
\headers{\ShortTitle}{\TheAuthors}
\title{{\TheTitle}}

\author{
Jordan T. Dworaczyk\thanks{School of Mathematical and Statistical Sciences, Arizona State University,
Tempe, AZ, United States (\email{jtdworac@asu.edu}, \email{malena.espanol@asu.edu})}
\and
Malena I. Espa\~nol$^*$}

\AddToHook{cmd/section/before}{\phantomsection}
\AddToHook{cmd/subsection/before}{\phantomsection}
\AddToHook{cmd/subsubsection/before}{\phantomsection}

\begin{document}
\maketitle

\begin{abstract}
Separable nonlinear least squares problems arise in inverse problems such as semi-blind image deblurring, where the image enters the forward model linearly, and the blur parameters enter nonlinearly. Even when the image is regularized, joint estimation of the image and blur parameters can favor a degenerate no-blur solution, in which the blur operator tends to the identity and the reconstruction merely denoises the data without removing the blur. We address this by regularizing
the blur parameters. Motivated by the determinant term arising from Gaussian marginalization, we construct a family of blur-parameter priors whose weight and auxiliary operators can be chosen independently of the image prior. The resulting log-determinant penalty is incorporated into regularized variable projection, which eliminates the image variables by minimization. Under periodic boundary conditions and assuming the matrices involved admit a joint spectral decomposition, we show that a sufficiently large determinant weight rules out the no-blur solution as a local minimizer. We also characterize the regularizer’s behavior as the auxiliary regularization scale tends to zero or infinity. We solve the reduced problem with a quasi-Newton method that evaluates the objective and its derivatives in the Fourier domain. In experiments with Gaussian, defocus, and motion blur, the proposed prior improves image and blur recovery over Gaussian marginalization across several settings, with determinant weights computed from the true image. A single-sample experiment illustrates how to automatically choose parameters from the data using a conditional generalized cross-validation criterion.
\end{abstract}



\section{Introduction}
\label{sec:intro}
We consider a semi-blind image deblurring problem, which we model by
\begin{equation}\label{eq: forward}
\bA(\by)\bx \approx \bb = \bb_{\mathrm{true}} + \boldsymbol{\epsilon}, \quad \text{with} \quad \bA(\by_{\mathrm{true}})\bx_{\mathrm{true}} = \bb_{\mathrm{true}},
\end{equation}
where the vector $\bb_{\mathrm{true}} \in \mathbb{R}^n$ represents unknown noise-free image data, while $\boldsymbol{\epsilon} \in \mathbb{R}^n$ represents additive measurement error, and $\bA(\by) \in \mathbb{R}^{n \times n}$ is a spatially invariant blurring operator with dependence on an unknown parameter vector $\by \in \calY\subset\mathbb{R}^r$. We assume $r \ll n$, that $\bA(\by)$ is severely ill-conditioned for $\by$ in a neighborhood of $\by_{\mathrm{true}}$, and that the mapping $\by \mapsto \bA(\by)$ is differentiable. Our objective is to estimate accurate approximations of both the image $\bx_{\mathrm{true}}$ and blur parameters $\by_{\mathrm{true}}$ from the noisy measurements $\bb$.

We formulate our objective as a \emph{regularized separable nonlinear least squares} (RSNLS) problem:
\begin{equation}\label{eq: rsnls}
\min_{\bx,\by}\phi(\bx,\by),\qquad\text{with}\quad \phi(\bx,\by) := \frac{1}{2}\|\bA(\by)\bx - \bb\|_2^2 + \frac{\lambda^2}{2}\|\bL\bx\|_2^2 + \frac{1}{2}\mathcal{R}(\by) ,
\end{equation}
where $\bL\in\mathbb{R}^{q\times n}$ is a regularization matrix (used to regularize $\bx$), $\lambda>0$ is the image regularization parameter, and $\calR(\by)$ is a blur-parameter penalty that includes its regularization weights.
We assume the null-space condition:
\begin{equation}\label{eq: null-space-condition}
\mathrm{null}(\bA(\by)) \cap \mathrm{null}(\bL) = \{\boldsymbol{0}\} \quad \text{for all } \by\in\calY.
\end{equation}
In the context of RSNLS, the term separable refers to the forward model $\bA(\by)\bx$ being linear in $\bx$ and nonlinear in $\by$.

For fixed blur parameters, Tikhonov regularization stabilizes the estimation of the image~\cite{deblurring}. However, joint image and blur estimation can favor a trivial solution: the estimated blurring matrix approaches the identity, while the estimated image retains the blur already present in the observations. This behavior is known as the \emph{no-blur degeneracy}, and image regularization alone does not necessarily prevent it.  This degeneracy has also been observed and analyzed in blind deconvolution~\cite{fergus2006removing,perrone2014total}. The authors in~\cite{levin2009understanding,levin2011efficient} distinguish joint estimation from blurring operator estimation after marginalizing over the image, showing how the latter can avoid the preference for no-blur solutions. In~\cite{wipf2014revisiting}, the authors also observe that marginalization prevents degeneracy and introduce a prior that yields a penalty term coupling the image, the blurring kernel, and the noise level.

Marginalization also underlies empirical Bayesian approaches to semi-blind deconvolution. In~\cite{mbakam2024marginal}, the blur parameters, noise level, and regularization parameters were estimated by maximizing the marginal likelihood obtained after integrating out the image. Then the reconstructed image is obtained by conditional maximum
a posteriori (MAP) estimation. The stochastic approximation method used addresses the generally intractable marginalization associated with non-Gaussian image priors.

We seek to construct a principled regularizer (or penalty term) $\mathcal{R}(\by)$ by examining the contribution of Gaussian marginalization.
Previous work has proposed Tikhonov-type and log-barrier penalty terms for $\by$~\cite{salzer2026variable}.

\medskip\noindent\textbf{Main contributions.}
\begin{enumerate}
    \item Motivated by Gaussian marginalization, we construct a family of priors for the blur parameters $\by$ and derive the corresponding regularizers $\calR(\by)$ for~\eqref{eq: rsnls}. The prior used to construct the penalty may be selected independently of the image reconstruction model. We establish conditions under which these priors are normalizable.

    \item For the case with periodic boundary conditions and operators diagonalized by the same discrete Fourier transform, we examine how the proposed penalty discourages the no-blur behavior and derive sufficient conditions under which the no-blur limit cannot be a local minimizer. We also characterize the penalty's behavior as the scale parameter in its construction approaches zero or infinity.

    \item We incorporate the proposed regularizer into the regularized variable projection method of \cite{salzer2026variable} and employ Fourier computations to evaluate the objective function and its derivatives. Experiments involving Gaussian, defocus, and motion blur assess image and blur recovery, sensitivity to initialization, and the effects of the regularization parameters. An additional experiment illustrates automatic selection of regularization parameters from the observed data.
\end{enumerate}

\medskip\noindent\textbf{Organization.}
Section~\ref{section: bayesian} develops the Bayesian model and the
proposed prior. Section~\ref{ssec:varpro-formulation} describes the VarPro
reduction of the RSNLS problem and analyzes the spectral behavior of the regularizer.
Section~\ref{sec: implementation details} describes the numerical
implementation, and Section~\ref{sec: numerics} presents the experiments.
Section~\ref{sec: Conclusions} concludes the paper.
The supplementary material gives derivative formulas, the Gauss-Newton derivation,
the automatic parameter-selection procedure, and additional reconstruction
metrics.

\section{Bayesian motivation and prior construction}
\label{section: bayesian}
In this section, we want to motivate and construct a blur-parameter regularizer $\mathcal R(\by)$. To do so, we first embed problem \eqref{eq: forward} within a Bayesian framework. We then use the determinant term arising from Gaussian marginalization to motivate a prior family for $\by$, which yields the corresponding regularizer. Consequently, we treat $\bx$ and $\by$ as random variables and apply Bayes' theorem to derive their joint \emph{posterior distribution} conditioned on the measured data $\bb$:
\begin{equation}\label{eq:posterior}
p\left(\bx,\by\, | \, \bb\right) = \frac{p\left(\bb\, |\, \bx,\by\right)p(\bx, \by)}{p(\bb)}.
\end{equation}
We consider the forward model~\eqref{eq: forward} under the assumption that the measurement error
is modeled by an isotropic Gaussian distribution $\bepsilon\sim\calN(\bzero, \eta^{-1}\bI_n)$. This gives the likelihood
$$ p(\bb\mid \bx, \by) \propto \exp\left(-\frac{\eta}{2}\|\bA(\by)\bx - \bb\|_2^2\right), $$
where $\eta>0$ denotes the precision of the random error. We hold $\eta$ fixed throughout this derivation. In the Supplementary Section~\ref{ssec:profiled-gcv-beta}, we propose an automatic parameter-selection procedure to estimate $\eta$ from $\bb$.

We assume that $\bx$ and $\by$ are a priori independent, so that the prior has the form $p(\bx, \by) = p(\bx) p(\by)$. To obtain the same image regularization term in ~\eqref{eq: rsnls}, we assume a Gaussian prior on $\bx$:
$$ p(\bx) \propto \exp\left(-\tfrac{\delta}{2}\|\bL\bx\|_2^2\right), $$
where $\delta>0$. When $\bL$ is rank deficient, this prior is improper and corresponds to an \emph{intrinsic} Gaussian distribution~\cite{rue2005gaussian}. In the deblurring setting, the rank-deficient case arises, for example, when a Gaussian Markov random field prior is chosen using a periodic discrete Laplacian stencil~\cite{rue2005gaussian,deblurring}.

For the unknown blur parameters, we leave the prior density $p(\by)$
unspecified for now. The semi-blind Bayesian model is thus summarized by
\begin{subequations}\label{eq: Bayesian-model}
    \begin{align}
    \bb\mid\bx, \by&\sim\calN(\bA(\by)\bx,\, \eta^{-1}\bI_n), \\
    p(\bx) &\propto \exp\!\left(-\tfrac{\delta}{2}\|\bL\bx\|_2^2\right), \\
    \by &\sim p(\by).
\end{align}
\end{subequations}
Combining this likelihood with the two priors gives the posterior distribution
$$ p(\bx, \by\mid \bb) \propto p(\bb\mid \bx, \by)p(\bx)p(\by) \propto \exp\left(-\frac{\eta}{2}\|\bA(\by)\bx - \bb\|_2^2 - \frac{\delta}{2}\|\bL\bx\|_2^2\right)p(\by). $$

Two common choices for $p(\by)$ are uniform and Gaussian priors. A uniform prior assigns a constant density over $\calY$, thereby assigning equal weight to all admissible blur parameters $\by$. Consequently, it imposes no penalty on the no-blur behavior. In contrast, an isotropic Gaussian prior centered at \(\mathbf y_0\in\calY\) has density
$$ p(\by)\propto\exp\!\left(-\frac{\delta_{\rm hp}}{2}\|\by-\by_0\|_2^2\right),\quad\delta_{\rm hp}>0, $$
and when increasing $\delta_{\rm hp}$, it concentrates the prior around $\by_0$, strengthening its influence and potentially overcoming this no-blur behavior~\cite{salzer2026variable}. However, in semi-blind deblurring, a reliable $\by_0$ close to $\by_{\rm true}$ may not be available.

\subsection{Gaussian marginalization and the determinant term}
\label{ssec:gaussian-structure}
We seek a prior that discourages the no-blur behavior without
assuming that the true blur is approximately known. To guide this
construction, we marginalize over the image in the linear-Gaussian
model presented above. The resulting determinant term depends on
the forward operator and motivates the prior family constructed below.

We begin with the conditional distribution $p(\bx\mid\by,\bb)$,
which is Gaussian under model~\eqref{eq: Bayesian-model}.
Its normalization will also provide the integral needed for marginalization.
To write this distribution explicitly, we define the precision-weighted
augmented matrix and data vector
$$ \bA_{\eta,\delta}(\by):=\begin{bmatrix}\sqrt\eta\,\bA(\by)\\[4pt]\sqrt\delta\,\bL\end{bmatrix},\qquad \bb_\eta:=\begin{bmatrix}\sqrt\eta\,\bb\\[4pt]\bzero\end{bmatrix}. $$
We also define
$$ \bQ(\by) := \eta\bA(\by)^\top\bA(\by)+\delta\bL^\top\bL = \bA_{\eta,\delta}(\by)^\top\bA_{\eta,\delta}(\by). $$
The null-space condition~\eqref{eq: null-space-condition} ensures
$\bQ(\by)$ is positive definite. For the calculations below, let
$$ \bA_{\eta,\delta}^\dagger(\by) = \big(\bA_{\eta,\delta}(\by)^\top\bA_{\eta,\delta}(\by)\big)^{-1} \bA_{\eta,\delta}(\by)^\top, \qquad \bmu_\bL(\by)=\bA_{\eta,\delta}^\dagger(\by)\bb_\eta, $$
and
$$ \bP_{\eta,\delta}^{\perp}(\by) = \bI-\bA_{\eta,\delta}(\by)\bA_{\eta,\delta}^\dagger(\by). $$
Here, $\bP^\perp_{\eta,\delta}(\by)$ is the orthogonal projector
onto the orthogonal complement of the column space of $\bA_{\eta,\delta}(\by)$.

Completing the square gives the identity
$$ \|\bA_{\eta,\delta}(\by)\bx-\bb_\eta\|_2^2 = \big(\bx-\bmu_\bL(\by)\big)^\top \bQ(\by) \big(\bx-\bmu_\bL(\by)\big) + \|\bP_{\eta,\delta}^{\perp}(\by)\bb_\eta\|_2^2, $$
which we can use to rewrite the joint posterior as
$$ p(\bx,\by\mid\bb)\propto \exp\!\left( -\frac{1}{2} \|\bP_{\eta,\delta}^\perp(\by)\bb_\eta\|_2^2 \right) \exp\left( -\frac12 \big(\bx-\bmu_\bL(\by)\big)^\top \bQ(\by) \big(\bx-\bmu_\bL(\by)\big) \right) p(\by). $$
Factors independent of $\bx$ cancel when the conditional posterior is normalized. Since $\bQ(\by)$ is positive definite, the Gaussian integral is
$$ \int_{\bbR^n} \exp\left( -\frac12 (\bx-\bmu_{\bL}(\by))^\top \bQ(\by) (\bx-\bmu_{\bL}(\by)) \right)\,d\bx = (2\pi)^{n/2}|\bQ(\by)|^{-1/2}. $$
The normalized conditional density is therefore
\begin{equation}\label{eq:conditional-Gaussian}
p(\bx\mid\by,\bb) = (2\pi)^{-n/2}|\bQ(\by)|^{1/2} \exp\left( -\frac12 (\bx-\bmu_{\bL}(\by))^\top \bQ(\by) (\bx-\bmu_{\bL}(\by)) \right),
\end{equation}
with mean $\bmu_{\bL}(\by)$ and covariance
$\bQ(\by)^{-1}$.

The same Gaussian integral eliminates $\bx$ from the joint posterior.
To normalize the resulting marginal density, define
$$ Z(\bb) := (2\pi)^{n/2} \int_{\calY} \exp\left( -\frac{1}{2}\|\bP_{\eta,\delta}^{\perp}(\by)\bb_\eta\|_2^2 \right) |\bQ(\by)|^{-1/2}p(\by)\,d\by, $$
and assume $0<Z(\bb)<\infty$. Integrating the normalized
joint posterior \eqref{eq:posterior} over $\bx$ gives
\begin{equation}\label{eq: marginal-posterior}
\begin{aligned}
p(\by\mid\bb) &= \int_{\bbR^n}p(\bx,\by\mid\bb)\,d\bx \\
&= \frac{(2\pi)^{n/2}}{Z(\bb)} |\bQ(\by)|^{-1/2} \exp\left( -\frac{1}{2}\|\bP_{\eta,\delta}^{\perp}(\by)\bb_\eta\|_2^2 \right)p(\by).
\end{aligned}
\end{equation}
To see how this determinant factor relates to the original joint problem,
combine~\eqref{eq:conditional-Gaussian} and~\eqref{eq: marginal-posterior}:
$$ p(\bx,\by\mid\bb) =p(\bx\mid\by,\bb)\,p(\by\mid\bb). $$
The prefactors in this product satisfy
$$ \underbrace{(2\pi)^{-n/2}|\bQ(\by)|^{1/2}}_{\text{conditional}} \underbrace{\frac{(2\pi)^{n/2}}{Z(\bb)}|\bQ(\by)|^{-1/2}}_{\text{marginal}} =\frac{1}{Z(\bb)}. $$
Thus, the determinant factors and the powers of $2\pi$ cancel when the conditional and marginal densities in~\eqref{eq:conditional-Gaussian} and~\eqref{eq: marginal-posterior} are recombined to recover the joint posterior. As discussed in the introduction, marginalizing over the image has been shown to mitigate the no-blur preference associated with joint MAP estimation.  Motivated by this behavior, we place a determinant factor directly in the blur prior so that its contribution remains in the joint posterior.

\subsection{The marginalization-based prior}
\label{ssec:marg-prior}

We define the \emph{marginalization-based prior family}
\begin{equation}\label{eq: prior-on-y}
p_\alpha(\by)\propto p_0(\by)|\tilde\bQ(\by)|^{-\alpha/2}, \qquad \by\in\calY,\quad \alpha\geq0.
\end{equation}
Here, $p_0$ is a base prior, such as a uniform or Gaussian density, which encodes prior information about $\by$. The determinant factor adjusts this density through the positive definite matrix $\tilde\bQ(\by)$, which we define by
$$ \tilde\bQ(\by) :=\tilde\eta\,\tilde\bA(\by)^\top\tilde\bA(\by) +\tilde\delta\,\tilde\bL^\top\tilde\bL = \tilde\bA_{\tilde\eta,\tilde\delta}(\by)^\top \tilde\bA_{\tilde\eta,\tilde\delta}(\by), $$
where
$$ \tilde\bA_{\tilde\eta,\tilde\delta}(\by) := \begin{bmatrix} \sqrt{\tilde\eta}\,\tilde\bA(\by)\\ \sqrt{\tilde\delta}\,\tilde\bL \end{bmatrix}. \qquad $$
We use $\tilde{\cdot}$ to denote parameters, matrices, and matrix operators that are \emph{auxiliary} to the linear-Gaussian model~\eqref{eq: Bayesian-model}. Notice that the auxiliary matrices $\tilde\bA(\by)\in\bbR^{m\times n}$ and $\tilde\bL\in\bbR^{\tilde q\times n}$ are weighted by fixed parameters $\tilde\eta>0$ and $\tilde\delta\geq0$, respectively.

When $\alpha=0$, we recover the base prior $p_0$, while $0<\alpha<1$ tempers how the determinant weights this density. The auxiliary operators and parameters can be chosen independently of the image model, allowing for choices where $\tilde\bQ(\by)\ne\bQ(\by)$. For example, $\tilde\bA(\by)$ can use periodic boundary conditions and the same PSF as a reconstruction operator $\bA(\by)$ with reflective boundaries. Together with a compatible periodic $\tilde\bL$, this makes $\tilde\bQ(\by)$ Fourier-diagonalizable for efficient determinant evaluation. In the analysis and numerical experiments below, we use periodic boundary conditions for both the reconstruction and auxiliary operators.

The following result establishes sufficient conditions under which this construction defines a proper prior.

\begin{proposition}[Properness of the marginalization-based prior]
Let $p_0$ be a proper density on $\calY$, and let
$\alpha\geq0$ and $\tilde\eta>0$ be fixed.
Suppose there exists a constant $c>0$ such that
$$ s_{\min}\big(\tilde\bA_{\tilde\eta,\tilde\delta}(\by)\big) \geq c \qquad\text{for all }\by\in\calY, $$
where $s_{\min}\big(\tilde\bA_{\tilde\eta,\tilde\delta}(\by)\big)$
denotes the smallest singular value of the auxiliary
augmented matrix $\tilde\bA_{\tilde\eta,\tilde\delta}(\by)$.
Then the normalizing integral
$$ Z_\alpha:=\int_{\calY}p_0(\by)|\tilde\bQ(\by)|^{-\alpha/2}\,d\by $$
satisfies $0<Z_\alpha<\infty$, and the normalizing prefactor of
the prior~\eqref{eq: prior-on-y} is $Z_\alpha^{-1}$.
\end{proposition}

\begin{proof}
Since
$\tilde\bQ(\by)=
\tilde\bA_{\tilde\eta,\tilde\delta}(\by)^\top
\tilde\bA_{\tilde\eta,\tilde\delta}(\by)$, we have
\(
|\tilde\bQ(\by)|\geq c^{2n}.
\)
Consequently,
$$ 0<Z_\alpha \leq c^{-\alpha n} \int_{\calY}p_0(\by)\,d\by = c^{-\alpha n} <\infty. $$
\end{proof}

\subsection{From MAP to regularized separable nonlinear least squares}

For the Bayesian model~\eqref{eq: Bayesian-model} with a blur prior
$p_\alpha$ of the form~\eqref{eq: prior-on-y}, taking the negative logarithm
of the joint posterior and omitting terms independent of $\bx$ and $\by$
gives the MAP problem
\begin{equation}\label{eq:MAP-estimation}
\min_{\bx,\by}\Bigl\{ \frac{\eta}{2}\|\bA(\by)\bx-\bb\|_2^2 +\frac{\delta}{2}\|\bL\bx\|_2^2 +\frac{\alpha}{2}\log\det\tilde\bQ(\by)-\log p_0(\by)\Bigr\}.
\end{equation}

Dividing this objective by the fixed positive precision $\eta$
preserves its minimizers and gives the regularization parameters
$$ \lambda^2=\frac{\delta}{\eta},\qquad \beta^2=\frac{\alpha}{\eta}. $$
Here, $\lambda$ controls image regularization and $\beta^2$
weights the determinant penalty.

We choose to rewrite the fixed scale $\tilde\eta$ within the log-determinant penalty of~\eqref{eq:MAP-estimation}, by writing $\tilde\bQ(\by)=\tilde\eta\,\tilde\bG_\gamma(\by)$, where
\begin{equation}\label{eq:auxiliary-gram}
\tilde\bG_\gamma(\by) :=\tilde\bA(\by)^\top\tilde\bA(\by) +\gamma\tilde\bL^\top\tilde\bL, \qquad \gamma:=\frac{\tilde\delta}{\tilde\eta}.
\end{equation}
The ratio $\gamma$ controls the relative contribution of the auxiliary regularizer. We therefore define the log-determinant and base-prior penalties by
$$ \mathcal R_\gamma(\by):=\log\det\tilde\bG_\gamma(\by), \qquad \mathcal{C}_0(\by):=-\eta^{-1}\log p_0(\by). $$
With these definitions, the scaled negative log-prior becomes
\begin{equation}\label{eq:constructed-regularizer}
-\eta^{-1}\log p_\alpha(\by) =\frac{\beta^2}{2}\mathcal R_\gamma(\by)+ \mathcal{C}_0(\by)+C,
\end{equation}
where $C:=\eta^{-1}\log Z_\alpha
+\frac{n\beta^2}{2}\log\tilde\eta$ is independent of $\by$. Here, the fixed scale $\tilde\eta$ contributes only an additive constant to $\log\det\tilde\bQ(\by)$. By omitting the constant $C$ in \eqref{eq:constructed-regularizer}, we identify the blur-parameter penalty
in~\eqref{eq: rsnls} as
$$ \calR(\by):=\beta^2\mathcal R_\gamma(\by)+2\mathcal{C}_0(\by). $$
The resulting deterministic RSNLS problem is
$\min_{\bx,\by}\phi(\bx,\by)$, where
\begin{equation}\label{eq:constructed-objective}
\phi(\bx,\by)=\frac12\|\bA(\by)\bx-\bb\|_2^2 +\frac{\lambda^2}{2}\|\bL\bx\|_2^2 +\frac{\beta^2}{2}\mathcal R_\gamma(\by)+\mathcal{C}_0(\by).
\end{equation}
With the regularization term $\mathcal{R}(\by)$ now specified, we turn to solving \eqref{eq:constructed-objective} via the variable projection method.

\section{Variable projection} \label{ssec:varpro-formulation}

The variable projection (VarPro) method was introduced by Golub and Pereyra~\cite{golub1973differentiation} for unregularized separable nonlinear least squares problems. It was recently extended in~\cite{salzer2026variable} to solve the RSNLS
problem~\eqref{eq: rsnls}, assuming that the regularization term $\mathcal R(\by)$ is differentiable. The resulting method is called \texttt{RGenVarPro}.

VarPro exploits the separable structure of the objective by minimizing over the linear variables $\bx$ for each fixed $\by$, leaving an optimization problem in the nonlinear parameters $\by$ alone. To express this reduction, we introduce the augmented matrix and data vector
$$ \bA_\lambda(\by):= \begin{bmatrix}\bA(\by)\\[4pt]\lambda\bL\end{bmatrix}, \qquad \bar\bb:= \begin{bmatrix}\bb\\[4pt]\bzero\end{bmatrix}, $$
and define
\begin{equation}\label{eq:deterministic-gram}
\begin{aligned}
\bG_\lambda(\by)&:=\bA_\lambda(\by)^\top\bA_\lambda(\by).
\end{aligned}
\end{equation}
The null-space condition~\eqref{eq: null-space-condition} guarantees that the matrix $\bA_\lambda(\by)$ has full column rank for every admissible $\by$. Then, for fixed $\by$, minimizing \eqref{eq:constructed-objective} with respect to $\bx$ yields the unique Tikhonov solution
\begin{equation}\label{eq:Tikhonov-solution}
\bx(\by)=\bigl(\bA(\by)^\top\bA(\by)+\lambda^2\bL^\top\bL\bigr)^{-1} \bA(\by)^\top\bb.
\end{equation}

Let
$$ \bP_\lambda^\perp(\by) :=\bI-\bA_\lambda(\by)\bG_\lambda(\by)^{-1} \bA_\lambda(\by)^\top $$
denote the orthogonal projector onto the orthogonal complement of
the column space of $\bA_\lambda(\by)$. Substituting $\bx(\by)$
into~\eqref{eq:constructed-objective} gives the reduced objective
$\bff_{\beta,\gamma}$ and the associated \emph{reduced problem}:
\begin{equation}\label{eq: reduced-problem}
\begin{aligned}
\bff_{\beta,\gamma}(\by)&:=\phi(\bx(\by),\by) =\varphi(\by) +\frac{\beta^2}{2}\mathcal R_\gamma(\by)+\mathcal{C}_0(\by), \\
\by^* &\in\argmin_{\by\in\calY}\bff_{\beta,\gamma}(\by),
\end{aligned}
\end{equation}
where 
$$ \varphi(\by) :=\frac12\|\bP_\lambda^\perp(\by)\bar\bb\|_2^2 $$
 is the \emph{reduced Tikhonov term}. Any minimizer $\by^*$ of the reduced problem yields a joint minimizer $\bigl(\bx(\by^*),\by^*\bigr)$ of the full objective $\phi(\bx,\by)$.

\begin{remark}
For given values of $\eta$ and $\delta$, and choosing $\tilde\bA=\bA$,
$\tilde\bL=\bL$, $\gamma=\lambda^2$, and $\beta^2=1/\eta$, the reduced objective $\phi(\bx(\by),\by)$ equals
$-\eta^{-1}\log p(\by\mid\bb)$ up to an additive constant independent
of $\by$, where $p(\by\mid\bb)$ is the marginal posterior
in~\eqref{eq: marginal-posterior} with prior $p(\by)=p_0(\by)$.
These choices define the Gaussian-marginalization benchmark used
in the numerical example presented in Section~\ref{ssec:reconstruction-experiments}.
\end{remark}

For further details on regularized variable projection and its
extensions, see~\cite{chen2018regularized,fu2019tikhonov,Espanol_2023,
chung2010efficient,gazzola2021regularization,pereyra2019imaging,
espanol2025local}. For constrained and nonsmooth formulations,
see~\cite{cornelio2014constrained,marmin2025varprox}.

\medskip

\subsection{Spectral analysis of the log-determinant regularizer}
\label{sec: analysis}

In this section, we analyze the effect of the log-determinant term on the reduced objective function produced by VarPro. To eliminate the effect of $\mathcal{C}_0(\by)$, we assume that it is constant and omit its contribution. We start by characterizing the no-blur reconstruction obtained by omitting $\mathcal{R}(\by)$, and its corresponding reduced objective value in a general setting. We then impose periodic boundary conditions and suitable spectral assumptions to compare this value with those associated with other blurs. This analysis reveals how the log-determinant term penalizes the no-blur limit and provides a sufficient condition on its weight to exclude that limit as a local minimizer. Finally, we examine how the auxiliary scale controls the penalty's dependence on the blur parameters.

We now characterize the no-blur reconstruction of the reduced problem by fixing $\lambda>0$, and setting $\beta=0$ in~\eqref{eq:constructed-objective}. Let $\overline{\calY}$ denote the closure of $\calY$. Suppose that there exists $\by^\times \in \overline{\calY}$, such that $\bA(\by)\to\bI$ as $\by\to\by^\times$. Under the null-space condition~\eqref{eq: null-space-condition}, the Tikhonov solution~\eqref{eq:Tikhonov-solution} is unique and satisfies
\begin{equation}\label{eq: degenerate-tikhonov-solution}
\bx(\by) \to \bx^\times := (\bI + \lambda^2\bL^\top\bL)^{-1}\bb, \qquad\text{as }\by\to\by^\times.
\end{equation}
The limiting reconstruction $\bx^\times$ is therefore the solution of a Tikhonov denoising problem, and the reduced objective~\eqref{eq: reduced-problem} converges to the no-blur objective value:
\begin{equation}\label{eq:degnerate-objective}
\begin{aligned}
\phi(\bx(\by),\by) &\to\min_{\bx}\left\{ \frac12\|\bx-\bb\|_2^2 +\frac{\lambda^2}{2}\|\bL\bx\|_2^2\right\} \\
&=\frac12\|\bx^\times-\bb\|_2^2 +\frac{\lambda^2}{2}\|\bL\bx^\times\|_2^2 \\
&=\frac12\bb^\top(\bb-\bx^\times).
\end{aligned}
\end{equation}
The last equality follows from
$(\bI+\lambda^2\bL^\top\bL)\bx^\times=\bb$.
Thus, in the no-blur limit, the reconstruction denoises the observations without correcting their blur.

We next compare how the reduced Tikhonov term $\varphi(\by)$ and the weighted log-determinant penalty change as the blur operator departs from the identity. For different choices of $\tilde\bG_\gamma(\by)$, this comparison yields sufficient conditions under which the no-blur limit is not a local minimizer of the reduced objective.

To simplify the analysis, we impose periodic boundary conditions and assume that the blurring
matrix $\bA(\by)$ and image regularization operator $\bL$ have
block-circulant-with-circulant-blocks (BCCB) structure. Both operators
are then simultaneously diagonalized by the unitary two-dimensional
discrete Fourier transform (DFT), denoted by $\bF$. That is, we have
\begin{equation}\label{eq:simultdiag2}
\bA(\by) = \bF^* \bC(\by) \bF, \qquad \bL = \bF^* \bS \bF,
\end{equation}
where
$$ \bC(\by)=\operatorname{diag}(\widehat c_0(\by),\ldots,\widehat c_{n-1}(\by)), \qquad \bS=\operatorname{diag}(\widehat s_0,\ldots,\widehat s_{n-1}). $$
The eigenvalues of each BCCB operator are obtained from the DFT of its first column; this is the standard Fourier diagonalization for periodic image deblurring~\cite[p.~42]{deblurring}. In particular, if $\ba_1(\by)$ denotes the first column of $\bA(\by)$ (equivalently, the vectorized zero-padded PSF arranged according to the BCCB convention), then $\widehat{\bc}(\by)=\sqrt{n}\,\bF\ba_1(\by)$, and the diagonal
entries of $\bS$ are obtained analogously from the first column of $\bL$.

\begin{remark}
The regularization operator $\bL$ need not be square or admit the representation $\bL=\bF^*\bS\bF$. It is sufficient that $\bL^\top\bL =\bF^*\operatorname{diag}(q_0,\ldots,q_{n-1})\bF$, with $q_k\geq0.$
In this case, $|\widehat{s}_k|^2$ in the Fourier formulas is replaced by $q_k$. For example, a rectangular operator
$\bL=[\bD_1^\top,\bD_2^\top]^\top$ satisfies this condition when $\bD_1$ and $\bD_2$ are periodic difference operators. The same observation applies to the auxiliary regularization operator $\tilde\bL$. We use this $\bL$ in the numerical experiments section.
\end{remark}

Define $\widehat{\bb}=\bF\bb$. Then, using the spectral decompositions of $\bA$ and $\bL$ given in \eqref{eq:simultdiag2}, we write the reduced Tikhonov term in the reduced problem \eqref{eq: reduced-problem} as
\begin{equation}\label{eq:phi0_vec}
\varphi(\by)= \frac12\sum_{k=0}^{n-1} \frac{\lambda^2|\widehat s_k|^2} {|\widehat c_k(\by)|^2+\lambda^2|\widehat s_k|^2} |\widehat b_k|^2.
\end{equation}
The Fourier representation~\eqref{eq:phi0_vec} explains why the RSNLS objective~\eqref{eq:constructed-objective} favors the no-blur limit when $\beta=0$.
For a nonnegative, unit-sum periodic
PSF, the triangle inequality gives
$|\widehat c_k(\by)|\leq1$ at every frequency.
Consequently, each summand in~\eqref{eq:phi0_vec}
is no smaller than its value at the identity. Then, we have
$$ \varphi(\by)\geq\varphi(\by^{\times}). $$
Thus, whenever the identity is an admissible operator or a limiting
member of the PSF family, $\varphi(\by)$ attains its minimum there
or approaches its infimum in the no-blur limit.

We next examine how the log-determinant term counteracts this preference.
First consider $\tilde\bG_\gamma(\by)=\bG_\lambda(\by)$, defined in \eqref{eq:deterministic-gram}, where $\tilde\bA(\by)=\bA(\by)$,
$\tilde\bL=\bL$, $\gamma=\lambda^2$. This is the case where the auxiliary and reconstruction operators coincide. Then, the reduced objective function is
$$ \bff_{\beta,\gamma}(\by) = \varphi(\by)+\frac{\beta^2}{2}\log\det\bG_\lambda(\by) = \varphi(\by)+\frac{\beta^2}{2} \sum_{k=0}^{n-1}\log g_k(\by), $$
where
\(
g_k(\by)
:=|\widehat c_k(\by)|^2+\lambda^2|\widehat s_k|^2
\)
are the Fourier eigenvalues of $\bG_\lambda$ in~\eqref{eq:deterministic-gram}.

We say that a blur exhibits \emph{spectral attenuation}
if
$|\widehat c_k(\by)|\leq1$ at every frequency, with
strict inequality at least once. Thus, the blur reduces
the magnitude of at least one Fourier component without
amplifying any others. The following propositions apply, under their stated assumptions,
to the general class of periodic PSFs satisfying this definition of
spectral attenuation.

\begin{proposition}[Spectral penalization of the no-blur limit]
\label{prop:spectral-penalization}
Assume $\lambda>0$. Let $\by^{\times}$
denote the no-blur limit, with
$|\widehat c_k(\by^{\times})|=1$ for every $k$.
Suppose $g_k(\by)>0$ at every frequency and
$$ |\widehat c_k(\by)|\leq1 \qquad\text{for every }k, $$
with strict inequality for at least one frequency. Then,
$$ \sum_{k=0}^{n-1}\log g_k(\by)< \sum_{k=0}^{n-1}\log g_k(\by^{\times}). $$
\end{proposition}

\begin{proof}
Since $t\mapsto\log(t+\lambda^2|\widehat s_k|^2)$ is strictly increasing on its domain, the
frequency-wise bounds
$|\widehat c_k(\by)|\leq1=|\widehat c_k(\by^{\times})|$ imply
$\log g_k(\by)\le \log g_k(\by^{\times})$ for each logarithmic summand, with strict
inequality in at least one summand. Summing over $k$ gives the result.
\end{proof}

Whereas Proposition~\ref{prop:spectral-penalization} shows that the proposed regularization penalizes no-blur behavior, the next proposition establishes a sufficient condition on the weight $\beta^2$ of the log-determinant term for excluding the no-blur solution as a local minimum, while permitting a more general auxiliary regularization operator. Consider the auxiliary choice
$\tilde\bA(\by)=\bA(\by)$, with fixed $\gamma\geq0$ and
an auxiliary operator diagonalized by the same Fourier transform:
$$ \tilde\bL=\bF^*\tilde\bS\bF, \qquad \tilde\bS=\operatorname{diag} (\widehat{\tilde s}_0,\ldots,\widehat{\tilde s}_{n-1}). $$
We allow $\tilde\bL$ to differ from $\bL$,
while retaining the reconstruction operator $\bL$ in $\varphi(\by)$.
Then,
$$ \bff_{\beta,\gamma}(\by)=\varphi(\by)+\frac{\beta^2}{2} \log\det\tilde\bG_\gamma(\by). $$

\begin{proposition}[Exclusion of the no-blur minimum]
\label{prop:no-blur-exclusion}
Assume the joint decomposition of~\eqref{eq:simultdiag2}, with fixed
$\lambda>0$, $\gamma\geq0$, and an auxiliary
operator $\tilde\bL$ diagonalized by $\bF$ as above. Assume
$\tilde\bA(\by)=\bA(\by)$. Let $\by^{\times}\in\overline{\calY}$
be a no-blur limit, and suppose the Fourier magnitudes extend continuously
to this point with $|\widehat c_k(\by^{\times})|=1$ for every $k$.
Suppose there are feasible parameters arbitrarily close to $\by^{\times}$
such that
$$ |\widehat c_k(\by)|\leq1\quad\text{for every }k, \qquad |\widehat c_j(\by)|<1\quad\text{for at least one }j. $$
Define $\ell_k:=\lambda^2|\widehat s_k|^2$ and
$d_k:=\gamma|\widehat{\tilde s}_k|^2$. Then $\by^{\times}$ is not a local minimizer of
$\bff_{\beta,\gamma}$ if
\begin{equation}\label{eq:no-blur-weight-bound}
\beta^2> \max_{0\leq k<n} \left[ (1+d_k) \frac{\ell_k|\widehat b_k|^2}{(1+\ell_k)^2} \right].
\end{equation}
\end{proposition}

\begin{proof}
We show that feasible blurs that are arbitrarily close to the identity have a smaller objective value. Set $a_k(\by)=|\widehat c_k(\by)|^2$, so that $a_k(\by^{\times})=1$ at every frequency.
The Fourier representation allows the reduced objective to be written as a sum over frequencies:
$$ \bff_{\beta,\gamma}(\by)=\sum_{k=0}^{n-1}h_k(a_k(\by)) \quad \mbox{ with } \quad h_k(a)=\frac12\left[ \frac{\ell_k|\widehat b_k|^2}{a+\ell_k} +\beta^2\log(a+d_k)\right]. $$
To determine whether decreasing $a$ from its no-blur
value of $1$ decreases $h_k(a)$, we differentiate it and obtain that
$$ h_k'(a)=\frac12\left[ -\frac{\ell_k|\widehat b_k|^2}{(a+\ell_k)^2} +\frac{\beta^2}{a+d_k}\right]. $$
The first term comes from the reduced Tikhonov term and the second from the log-determinant term. At $a=1$,
condition~\eqref{eq:no-blur-weight-bound} ensures
$$ h_k'(1)=\frac12\left[ -\frac{\ell_k|\widehat b_k|^2}{(1+\ell_k)^2} +\frac{\beta^2}{1+d_k}\right]>0 \qquad\text{for every }k. $$
At the identity, the log-determinant contribution therefore outweighs the contribution of $\varphi(\by)$ in every spectral derivative. We now use continuity to show that this balance persists for nearby blurs. Spectral attenuation will then decrease at least one frequency contribution without increasing the others, giving a smaller total objective. Since $h_k'(1)>0$ and $h_k'(a)$ is continuous, for each frequency $k$ there exists $\epsilon_k>0$ such that $h_k'(a)>0$ whenever $|a-1|<\epsilon_k$. There are finitely many frequencies, so take $\epsilon=\min_k\epsilon_k>0$. Then, for every $k$, $h_k(a)$ is strictly increasing on $(1-\epsilon,1+\epsilon)$.

Since $a_k(\by)\to1$ as $\by\to\by^{\times}$, there exists a neighborhood of $\by^{\times}$ in
which 
$$ |a_k(\by)-1|<\epsilon\qquad\text{for every }k. $$
Choose a feasible $\by$ sufficiently close to $\by^{\times}$ that exhibits spectral attenuation. Then, $a_k(\by)\leq1$ for every $k$, with strict inequality for at least one frequency. Because each $h_k$ is strictly increasing on $(1-\epsilon,1+\epsilon)$,
$$ h_k(a_k(\by))\leq h_k(1) \qquad\text{for every }k, $$
with at least one strict inequality. Consequently,
$$ \bff_{\beta,\gamma}(\by) =\sum_{k=0}^{n-1}h_k(a_k(\by)) < \sum_{k=0}^{n-1}h_k(1) =\bff_{\beta,\gamma}(\by^{\times}). $$
By assumption, feasible parameters exhibiting spectral attenuation can be chosen arbitrarily close to $\by^{\times}$. The preceding inequality shows that these parameters have smaller objective values. Hence, $\by^{\times}$ is not a local minimizer.
\end{proof}

\begin{remark}
The Gaussian PSFs in Section~\ref{sec: numerics} approach the identity as their widths tend to zero. The smoothed motion and defocus families retain residual blur at fixed smoothing parameters, so Proposition~\ref{prop:no-blur-exclusion} does not directly cover their low-blur limits.
\end{remark}

The exclusion threshold~\eqref{eq:no-blur-weight-bound} shows how the auxiliary scale $\gamma$ enters the condition for excluding the no-blur limit. We now examine how the log-determinant penalty behaves as $\gamma$ tends to zero or infinity. Using $\tilde\bG_\gamma(\by)$ from~\eqref{eq:auxiliary-gram}, we define
$$ \mathcal R_{\beta}(\by;\gamma):=\frac{\beta^2}{2}\log\det\tilde\bG_\gamma(\by). $$

\begin{proposition}[Asymptotic regimes of the auxiliary log-determinant]
\label{prop:delta_limits}
Fix $\by$, $\beta^2\geq0$, and the
auxiliary operators $\tilde\bA(\by)$ and $\tilde\bL$.
If $\tilde\bL^\top\tilde\bL$ is positive definite, then
\begin{equation}\label{eq:delta-large-regime}
\mathcal R_{\beta}(\by;\gamma) = \frac{\beta^2 n}{2}\log\gamma +\frac{\beta^2}{2}\log\det(\tilde\bL^\top\tilde\bL) +\mathcal O(\gamma^{-1}), \qquad \gamma\to\infty.
\end{equation}
If $\tilde\bA(\by)^\top\tilde\bA(\by)$ is positive definite, then
\begin{equation}\label{eq:delta-small-regime}
\mathcal R_{\beta}(\by;\gamma) = \frac{\beta^2}{2} \log\det(\tilde\bA(\by)^\top\tilde\bA(\by)) +\mathcal O(\gamma), \qquad \gamma\to0.
\end{equation}
\end{proposition}

\begin{proof}
For the large-scale limit, factoring
$\gamma\,\tilde\bL^\top\tilde\bL$ and taking
log-determinants gives
$$ \log\det\tilde\bG_\gamma(\by) =n\log\gamma +\log\det(\tilde\bL^\top\tilde\bL) +\log\det\left[\bI+\frac1\gamma (\tilde\bL^\top\tilde\bL)^{-1}\tilde\bA(\by)^\top\tilde\bA(\by) \right]. $$

\begingroup\emergencystretch=1em
Using $\log\det(\bI+\epsilon\bM)=\mathcal O(\epsilon)$
for fixed matrix $\bM$, and multiplying by
$\beta^2/2$, yields~\eqref{eq:delta-large-regime}. Similarly, factoring
$\tilde\bA(\by)^\top\tilde\bA(\by)$ gives\par\endgroup
$$ \log\det\tilde\bG_\gamma(\by) =\log\det(\tilde\bA(\by)^\top\tilde\bA(\by))+\log\det\left[ \bI+\gamma \big(\tilde\bA(\by)^\top\tilde\bA(\by)\big)^{-1} \tilde\bL^\top\tilde\bL \right]. $$

The same expansion yields~\eqref{eq:delta-small-regime}.
\end{proof}

Under the assumptions of
Proposition~\ref{prop:delta_limits}, the leading terms
in~\eqref{eq:delta-large-regime} are independent of
$\by$, whereas~\eqref{eq:delta-small-regime} shows
that the small-scale limit is governed by the spectrum
of the auxiliary forward operator. Thus, a large
$\gamma$ can weaken the regularizer's dependence
on the blur parameters, while a small $\gamma$
can make it sensitive to small eigenvalues of
$\tilde\bA(\by)^\top\tilde\bA(\by)$. Such sensitivity
may produce large gradients and curvature, as
the derivative formulas in
Supplementary Section~\ref{app:logdet-derivatives} illustrate. These limiting regimes motivate testing intermediate
auxiliary scales. For $\gamma>0$, choosing $\tilde\bL=\bI$ gives
$\tilde\bG_\gamma$ a uniform eigenvalue lower bound of
$\gamma$, limiting sensitivity to small spectral
values while retaining dependence on the blur.
Section~\ref{ssec:auxiliary-prior-sensitivity}
examines this choice with $\gamma=1/\eta$ and
evaluates the effect of varying $\gamma$.

\section{Quasi-Newton method and implementation details}
\label{sec: implementation details}
We solve the reduced problem~\eqref{eq: reduced-problem} using
a quasi-Newton method based on~\cite{salzer2026variable},
with Hessian damping and an Armijo line search. Recall that, for fixed
regularization parameters, the reduced objective is
\begin{equation}\label{eq:implementation-objective}
\bff_{\beta,\gamma}(\by) =\varphi(\by) +\frac{\beta^2}{2}\mathcal R_\gamma(\by) +\mathcal C_0(\by).
\end{equation}
We first describe the quasi-Newton iteration and then explain how its
components are evaluated using the Fourier representation.

The iterates are defined by
$$ \by_{k+1}=\by_k+t_k\bp_k, $$
where $t_k>0$ is the step length and $\bp_k$ is the search
direction, computed by solving
\begin{equation}\label{eq:damped-step}
(\bH(\by_k)+\tau_k\bI)\bp_k =-\nabla\bff_{\beta,\gamma}(\by_k).
\end{equation}
Here, $\tau_k\geq0$ is a damping parameter, and the Hessian
approximation is
\begin{equation}\label{eq:step-matrix}
\bH(\by) =\bH_{\mathrm{GN}}(\by) +\beta^2\bH_{\log}(\by) +\bH_{\mathrm{prior}}(\by),
\end{equation}
where $\bH_{\mathrm{GN}}$ is a Gauss-Newton approximation
of $\nabla^2\varphi(\by)$, and
$$ \bH_{\log}(\by) =\frac12\nabla^2\mathcal R_\gamma(\by), \qquad \bH_{\mathrm{prior}}(\by) =\nabla^2\mathcal C_0(\by). $$
For a uniform base prior,
$\bH_{\mathrm{prior}}=\boldsymbol{0}$.
For an isotropic Gaussian base prior with precision
$\delta_{\rm hp}\bI$, we have
$\bH_{\mathrm{prior}}=w_{\rm hp}\bI$, where
$w_{\rm hp}=\delta_{\rm hp}/\eta$.

The damping parameter $\tau_k$ is chosen to make the
coefficient matrix in~\eqref{eq:damped-step} positive definite.
The determinant Hessian can be indefinite even when
$\tilde\bG_\gamma$ is positive definite, so damping may
be necessary. We first attempt Cholesky factorization
with $\tau_k=0$. If the factorization fails, we increase
$\tau_k$ by successive increments
$\tau_{\mathrm{base}},\kappa\tau_{\mathrm{base}},
\kappa^2\tau_{\mathrm{base}},\ldots$,
where $\tau_{\mathrm{base}}>0$ and $\kappa=10$.
If no factorization succeeds within the prescribed number
of attempts, we use the negative-gradient direction
$\bp_k=-\nabla\bff_{\beta,\gamma}(\by_k)$.
For $\nabla\bff_{\beta,\gamma}(\by_k)\ne\bzero$,
either choice gives a descent
direction~\cite[Chs.~3 and~6]{nocedal2006numerical}.

Given $\bp_k$, Armijo backtracking starts from an initial
step length $t_{\mathrm{init}}>0$ and reduces it by a factor
$\zeta\in(0,1)$ until the trial parameter is feasible and
\begin{equation}\label{eq:armijo-v14}
\bff_{\beta,\gamma}(\by_k+t_k\bp_k) \leq \bff_{\beta,\gamma}(\by_k) +c_1t_k[\nabla\bff_{\beta,\gamma}(\by_k)]^\top\bp_k, \qquad 0<c_1<1.
\end{equation}
The iteration stops when the gradient norm or the accepted
step norm falls below its prescribed tolerance, when
backtracking fails, or when the maximum number of iterations
is reached. Algorithm~\ref{alg:HybridGN} summarizes the procedure.

\begin{algorithm}
\caption{Quasi-Newton method with damping and line search}
\label{alg:HybridGN}
\begin{algorithmic}[1]
\Require Feasible $\by_0$, fixed regularization parameters,
tolerances $\epsilon_g,\epsilon_s$, iteration budget $K$,
and damping and line-search parameters
\State $\by\gets\by_0$
\State Set termination flag to iteration limit
\For{$k=0,\ldots,K-1$}
    \State Evaluate $\bff_{\beta,\gamma}(\by)$ and
    $\bg\gets\nabla\bff_{\beta,\gamma}(\by)$
    \If{$\|\bg\|_2\leq\epsilon_g$}
        \State Set gradient flag
        \State \textbf{break}
    \EndIf
    \State Form  $\bH_k$ using~\eqref{eq:step-matrix}
    \State Compute $\bp_k$ using~\eqref{eq:damped-step},
    with damping and the negative-gradient fallback
    described above
    \State Apply backtracking to find a feasible $t_k$
    satisfying~\eqref{eq:armijo-v14}
    \If{backtracking fails}
        \State Set step flag (zero step)
        \State \textbf{break}
    \EndIf
    \State $\by^+\gets\by+t_k\bp_k$
    \State $s_k\gets\|\by^+-\by\|_2$
    \State $\by\gets\by^+$
    \If{$s_k\leq\epsilon_s$}
        \State Set step flag
        \State \textbf{break}
    \EndIf
\EndFor
\State Evaluate the objective and gradient at the final $\by$
\State \Return $\by$, objective value, gradient norm,
and termination flag
\end{algorithmic}
\end{algorithm}

\subsection{Fourier evaluation of the objective and derivatives}
\label{ssec:fft-evaluation}

For this evaluation, we use the Fourier diagonalization
in~\eqref{eq:simultdiag2} and specialize to $\tilde\bA=\bA$,
with $\tilde\bL$ diagonalized by the same transform.
The symbols $\widehat{\tilde s}_k$ denote the eigenvalues of $\tilde\bL$.
For the Hessian formulas, we assume twice continuously differentiable dependence of $\bA$ on $\by$.
We retain the data convention $\widehat\bb=\bF\bb$
from Section~\ref{sec: analysis}. In \textsc{Matlab}, these coefficients
are obtained by multiplying the vectorized output of \texttt{fft2}
by $1/\sqrt n$.

For fixed $\lambda$, each component of the gradient of $\varphi(\by)$ is
\begin{equation}\label{eq:fft-reduced-gradient}
\partial_{y_i} \varphi(\by) =-\sum_{k=0}^{n-1} \frac{\lambda^2|\widehat s_k|^2|\widehat b_k|^2 \operatorname{Re}(\overline{\widehat c_k}\, \partial_{y_i}\widehat c_k)} {(|\widehat c_k|^2+\lambda^2|\widehat s_k|^2)^2}
\end{equation}
and the corresponding Gauss–Newton approximation is defined by
\begin{equation}\label{eq:fft-gn-matrix}
(\bH_{\mathrm{GN}})_{i\ell} =\sum_{k=0}^{n-1} \frac{\lambda^2|\widehat s_k|^2|\widehat b_k|^2 \operatorname{Re}(\overline{\widehat c_k}\, \partial_{y_i}\widehat c_k) \operatorname{Re}(\overline{\widehat c_k}\, \partial_{y_\ell}\widehat c_k)} {(|\widehat c_k|^2+\lambda^2|\widehat s_k|^2)^3}.
\end{equation}
This approximation uses the Fourier residual vector~\eqref{eq:compressed-residual} and its Jacobian~\eqref{eq:compressed-jacobian}, derived in Supplementary Section~\ref{app:varpro-gn}.

The same Fourier representation gives the log-determinant contribution
in~\eqref{eq:implementation-objective}. The auxiliary matrix
$\tilde\bG_\gamma$ has eigenvalues
$$ \tilde g_k(\by) =|\widehat c_k(\by)|^2+\gamma|\widehat{\tilde s}_k|^2. $$
Its one-half log-determinant, gradient, and Hessian are evaluated as
$$
\begin{aligned}
\frac12\log\det\left( \bA^\top\bA+\gamma\tilde\bL^\top\tilde\bL\right) &=\frac12\sum_{k=0}^{n-1}\log\tilde g_k, \\
\partial_{y_i}\left(\frac12\log\det\tilde\bG_\gamma\right) &=\frac12\sum_{k=0}^{n-1} \frac{\partial_{y_i}\tilde g_k}{\tilde g_k}, \\
(\bH_{\log})_{i\ell} &=\frac12\sum_{k=0}^{n-1}\left[ \frac{\partial^2_{y_i y_\ell}\tilde g_k}{\tilde g_k} -\frac{(\partial_{y_i}\tilde g_k) (\partial_{y_\ell}\tilde g_k)}{\tilde g_k^2}\right].
\end{aligned}
$$
The required eigenvalue derivatives follow from
$$
\begin{aligned}
\partial_{y_i}\tilde g_k &=2\operatorname{Re}(\overline{\widehat c_k}\, \partial_{y_i}\widehat c_k), \\
\partial^2_{y_i y_\ell}\tilde g_k &=2\operatorname{Re}\left( \overline{\partial_{y_\ell}\widehat c_k}\, \partial_{y_i}\widehat c_k +\overline{\widehat c_k}\, \partial^2_{y_i y_\ell}\widehat c_k\right).
\end{aligned}
$$
These expressions specialize the trace identities in Supplementary Section~\ref{app:logdet-derivatives}.
The Fourier symbols and their derivatives are obtained by applying FFTs to the PSF and its parameter derivatives. For a fixed small parameter dimension $r$, evaluation requires $\mathcal O(n\log n)$ arithmetic and $\mathcal O(n)$ storage.

When such a spectral decomposition is unavailable, the inexact RGenVarPro (\texttt{iRGenVarPro}) method introduced in~\cite{salzer2026variable} could be combined with approaches for approximating the log-determinant term, such as the majorization-minimization scheme proposed in~\cite{buser2026majorization}.

\section{Numerical experiments}
\label{sec: numerics}
In this section, we present numerical experiments to study the effect of the proposed log-determinant regularizer in semi-blind image deblurring. We solve the corresponding RSNLS problem using \texttt{RGenVarPro} with the quasi-Newton method described in Section~\ref{sec: implementation details}. We first examine the objective landscapes and then compare image and blur recovery for different choices of auxiliary matrices and initial guesses. Next, we study the effect of the auxiliary scale. Finally, we examine the performance of automatic parameter selection. Our implementations of the methods described in this paper will be available at \url{https://github.com/jordan-dworaczyk/MargRegVarPro}.

All experiments assume periodic boundary conditions and spatially invariant blur. We add Gaussian noise at the stated relative noise level, and the true image intensities lie in $[0,1]$. In the first three experiments, we prescribe the noise precision and image regularization parameter. In the final experiment, we estimate these parameters from the observations using the procedure in Supplementary Section~\ref{ssec:profiled-gcv-beta}. We report the prior exponent $\alpha$, which determines the weight $\beta^2=\alpha/\eta$.

To define the blurring operators, we consider three classes of parameterized point spread functions (PSFs), corresponding to Gaussian, motion, and defocus blur. Kernels use centered row and column coordinates $(s,t)$; in each definition, $Z$ denotes the parameter-dependent normalization factor making the discrete kernel sum to one. The smooth parameterizations permit derivative-based optimization.

\medskip

\paragraph{Gaussian PSF}
The anisotropic Gaussian blur is modeled using the kernel
\begin{equation}\label{eq:GaussianPSF}
g_{\sigma_1,\sigma_2}(s,t) = \frac{1}{Z(\sigma_1,\sigma_2)} \exp\!\left(-\frac{s^2}{\sigma_1^2} - \frac{t^2}{\sigma_2^2}\right),
\end{equation}
where $(\sigma_1,\sigma_2) > 0$ are the blur parameters. For this PSF, $\by=(\sigma_1,\sigma_2)$.

\medskip

\paragraph{Motion PSF}
The motion blur PSF corresponds to a uniform line segment of
length $L>0$ and orientation $\theta\in\mathbb{R}$. We define a smooth motion blur by letting
$$ u(s,t;\theta) = t\cos\theta - s\sin\theta\qquad\text{and} \qquad v(s,t;\theta) = t\sin\theta + s\cos\theta $$
denote the coordinates along and orthogonal to the direction of motion, respectively. Thus, $\theta=0$ gives horizontal motion. The motion PSF is then written as
$$
\begin{aligned}
g_{L,\theta}(s,t) &= \frac{1}{Z(L,\theta)} \exp\!\left( -\frac{v(s,t;\theta)^2}{2\sigma_{\perp}^2} \right) \\
&\quad \times \frac{1}{2} \left[ \tanh\!\left( \frac{u(s,t;\theta)+L/2}{\tau} \right) - \tanh\!\left( \frac{u(s,t;\theta)-L/2}{\tau} \right) \right],
\end{aligned}
$$
where $\sigma_{\perp}=0.30$ controls the transverse width and $\tau=0.15$ smooths the endpoints. For this PSF, $\by=(L,\theta)$.

\medskip

\paragraph{Defocus PSF}
The defocus blur is modeled as a smoothed disk with radius $R > 0$.
Let $r = \sqrt{s^2 + t^2}$ denote the radial coordinate. The PSF is given by
\begin{equation}\label{eq:defocus-psf}
g_R(s,t) = \frac{1}{Z(R)} \frac{1}{1 + \exp\!\left(\frac{r - R}{\varepsilon_{\mathrm{edge}}}\right)},
\end{equation}
where $\varepsilon_{\mathrm{edge}}=0.5$ is a fixed smoothing parameter controlling the width of the transition at the disk boundary. In this PSF, the parameter is a scalar, so $\by=R.$

\subsection{Objective landscape comparisons}
We examine how different combinations of regularization terms for $\by$ change the landscape of the reduced objective.
We consider the Gaussian PSF defined in~\eqref{eq:GaussianPSF} and image regularization based on periodic second differences
in the two coordinate directions. Let $\bD_1$ and $\bD_2$ denote these operators. We choose $\bL$ so that
$$ \bL^\top\bL=\bD_1^\top\bD_1+\bD_2^\top\bD_2, $$
which gives
$\|\bL\bx\|_2^2=\|\bD_1\bx\|_2^2+\|\bD_2\bx\|_2^2$.

For regularization of $\by$, we consider an $\ell_2$-norm regularization centered at an inaccurate reference parameter
$\by_c$, with weight $w_{\mathrm{hp}}$, and a log-determinant regularizer with matched auxiliary matrices,
$\tilde\bG_\gamma=\bG_\lambda$. The reduced objective is therefore
$$ \bff_{\beta,\gamma}(\by) =\varphi(\by) +\frac{\beta^2}{2}\log\det\bG_\lambda(\by) +\frac{w_{\mathrm{hp}}}{2}\|\by-\by_c\|_2^2. $$
We compare the objective landscapes for all six combinations of $\beta^2\in\{0,1/\eta\}$ and $w_{\mathrm{hp}}\in\{0,10^2/\eta,10^3/\eta\}$. The experimental setup is summarized in Table~\ref{tab:gauss-exp-params}.

\begin{table}[t]
    \centering
    \caption{Setup for the Gaussian PSF objective-landscape experiment.}
    \label{tab:gauss-exp-params}
    \begin{tabular}{@{}ll@{}}
        \toprule
        Quantity & Value \\
        \midrule
        Test image & \emph{Cameraman}, $128 \times 128$ \\
        PSF model & Anisotropic Gaussian \\
        True PSF parameter & $\by_{\mathrm{true}}=(3,5)$ \\
        Relative noise level & $1\%$ \\
        Image penalty & Periodic second differences in both directions \\
        Image regularization & $\lambda^2=100/\eta$ \\
        Parameter-prior mean & $\by_c=(4,4)$ \\
        Parameter-prior weight & $w_{\mathrm{hp}}\in\{0,10^2/\eta,10^3/\eta\}$ \\
        Determinant weight & $\beta^2\in\{0,1/\eta\}$ \\
        Landscape discretization & $80\times80$ uniform grid on $[0.01,7]^2$ \\
        \bottomrule
    \end{tabular}
\end{table}

\begin{figure}[t]
\centering
\includegraphics[width=\linewidth]{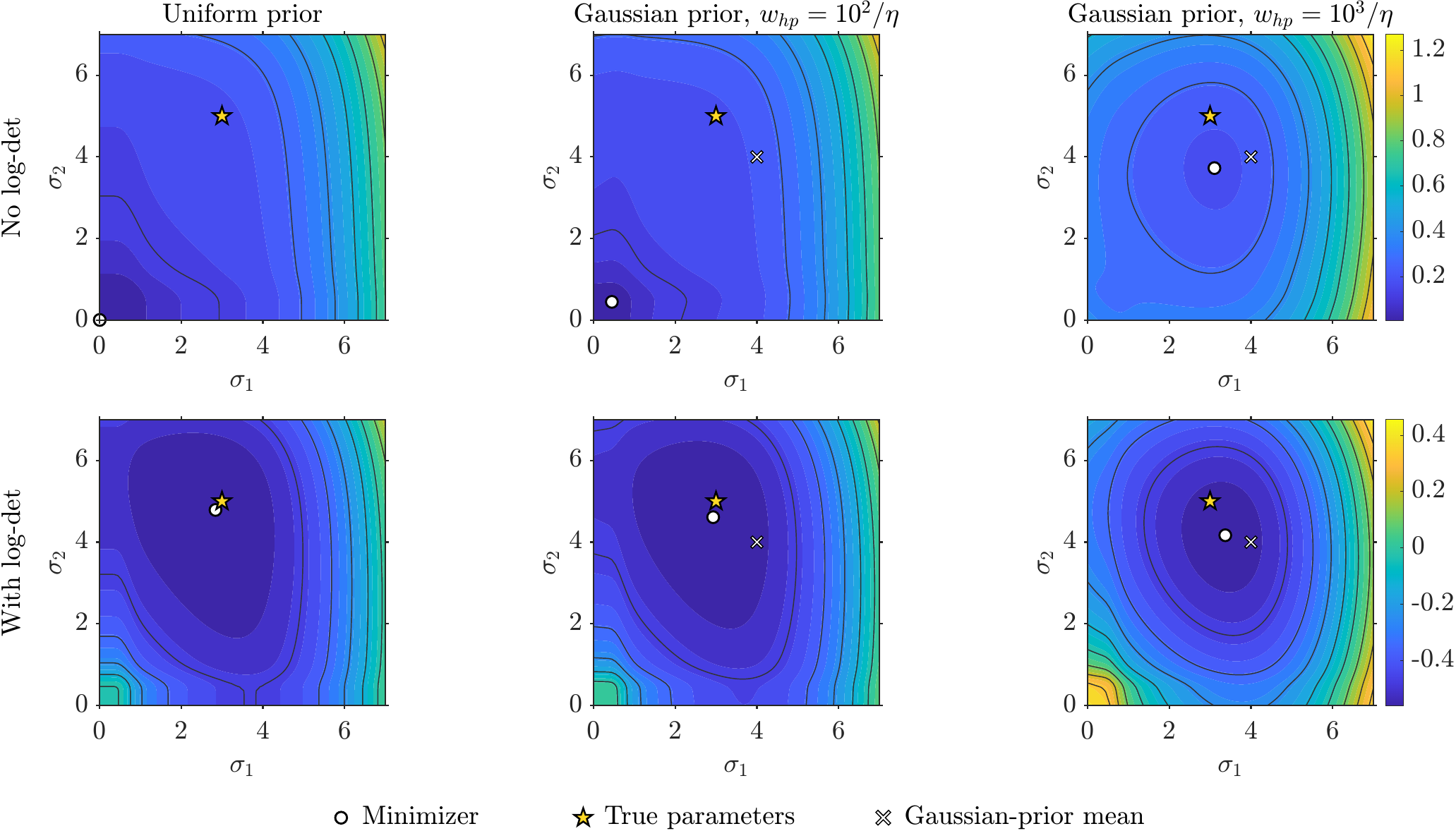}
\caption{\emph{Effect of prior choice on the objective landscape for blur parameter estimation.} The rows labeled ``No log-det'' and ``With log-det'' correspond to $\beta^2=0$ and $\beta^2=1/\eta$, respectively, with the three base priors indicated by the column headings. Circles mark grid minimizers, stars mark the true blur parameters, and outlined white crosses mark the Gaussian-prior mean. Raw objective values share a common color scale within each row.}
\label{fig:Gaussian-vs-uniform}
\end{figure}

\medskip
\noindent\textbf{Results.} Figure~\ref{fig:Gaussian-vs-uniform} shows the objective landscapes
for the six combinations of regularization weights. For $\beta^2=0$ and $w_{\mathrm{hp}}=0$, the minimum over the
plotting grid occurs at $(0.01,0.01)$, close to the no-blur limit. Adding the quadratic penalty with $w_{\mathrm{hp}}=10^2/\eta$ moves the grid minimizer to approximately $(0.45,0.45)$. Increasing the weight to $10^3/\eta$ moves it to approximately $(3.11,3.73)$. Thus, a sufficiently strong quadratic penalty moves the minimizer away from the no-blur limit, but its location depends on the reference parameter $\by_c=(4,4)$. For $\beta^2=1/\eta$, the approximate grid minimizers are $(2.84,4.79)$, $(2.93,4.61)$, and $(3.37,4.17)$ for $w_{\mathrm{hp}}=0$, $10^2/\eta$, and $10^3/\eta$, respectively. Without the quadratic penalty, the log-determinant term moves the grid minimizer close to the true parameter $\by_{\mathrm{true}}=(3,5)$. This behavior is consistent with Proposition~\ref{prop:spectral-penalization}, which shows that the log-determinant term penalizes the no-blur limit relative to blurs exhibiting spectral attenuation. As the quadratic weight increases, the grid minimizer moves toward $\by_c$ and farther from $\by_{\mathrm{true}}$ in this example. Thus, although the log-determinant term counteracts the preference for no blur, the quadratic penalty can still shift the minimizer toward an inaccurate reference parameter.

\subsection{Image reconstruction experiments}
\label{ssec:reconstruction-experiments}
\begin{figure}[t]
\centering
\includegraphics[width=\linewidth]{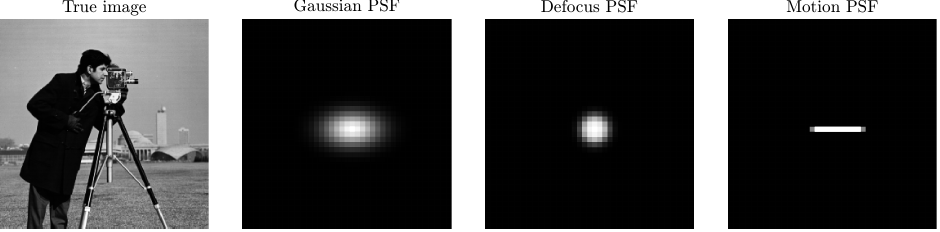}
\caption{\emph{True image and PSFs used in the reconstruction experiments.}
The left panel shows the true Cameraman image. The remaining
panels show the true Gaussian, defocus, and motion PSFs used
to generate the corresponding blurred and noisy data.
Each PSF is shown on a cropped support and normalized
independently for display.}
\label{fig:ground-truth-reference}
\end{figure}

\begin{figure}[t]
\centering
\includegraphics[width=\linewidth]{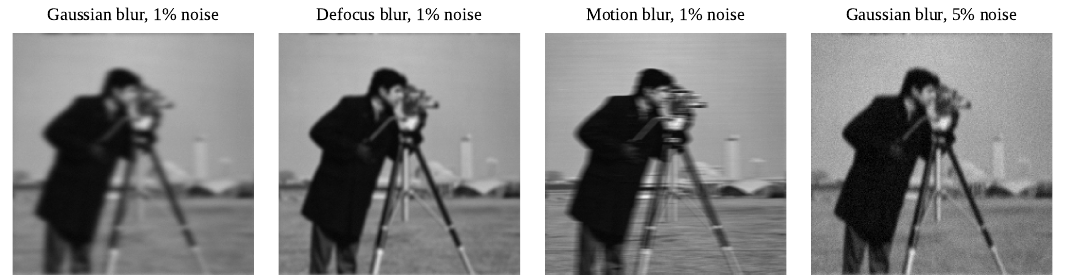}
\caption{\emph{Blurred and noisy data.} The first three panels show the observed data for the Gaussian, defocus, and motion blurs of Figure~\ref{fig:ground-truth-reference} with $1\%$ relative noise (reconstruction experiments). The last panel shows the Gaussian blur with $5\%$ relative noise used in the automatic parameter-selection experiment. All images use the intensity scale $[0,1]$.}
\label{fig:blurred-data}
\end{figure}

\raggedbottom
We now compare image and blur reconstructions for three choices of regularization: image regularization alone, the log-determinant regularization with matched auxiliary matrices (i.e., $\tilde\bG_\gamma=\bG_\lambda$), and the log-determinant regularization with $\tilde\bA=\bA$ and $\tilde\bL=\bI$. We consider Gaussian, defocus, and motion PSFs, with true Gaussian widthts $(3,5)$, radius $2.5$, and length-angle pair $(10,0)$, respectively. We use the $256\times256$ \emph{Cameraman} image with $1\%$ relative noise. The true image and PSFs are shown in Figure~\ref{fig:ground-truth-reference}, and the resulting blurred and noisy data are shown in Figure~\ref{fig:blurred-data}. For image regularization, we use the operator $\bL$ described in the previous experiment, based on periodic second differences in both coordinate directions, with $\lambda^2=100/\eta$. We assume that $\mathcal C_0(\by)$ in~\eqref{eq:implementation-objective} is constant and omit its contribution. We denote the number of image pixels by $N=256^2$.

\begin{figure}[!htbp]
\centering
\includegraphics[width=\linewidth]{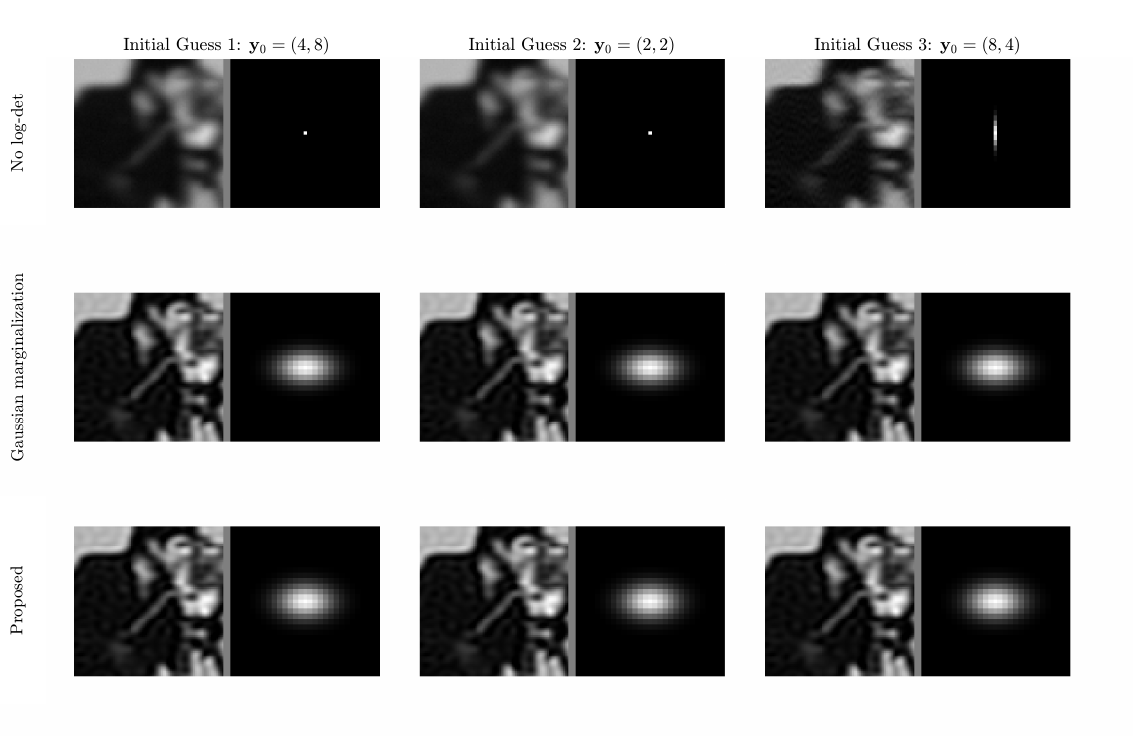}
\caption{\emph{Gaussian blur: paired zoomed reconstructions and recovered PSFs from all initializations.}
Each cell pairs a zoomed reconstruction (left) with its recovered PSF (right). Images use the scale $[0,1]$ and each PSF is peak-normalized for display. Rows are labeled ``No log-det'', ``Gaussian marginalization'', and ``Proposed''; columns show the three initial guesses in Table~\ref{tab:reconstruction-gaussian}.}
\label{fig:composite-recon-psf-gaussian}
\par\medskip
\begin{minipage}{\linewidth}
\centering
\normalsize
\resulttablecaption{Recovery for Gaussian blur from all three initializations (proposed $\alpha=0.5$).}
\label{tab:reconstruction-gaussian}
\setlength{\tabcolsep}{4.5pt}
\renewcommand{\arraystretch}{1.05}
\begin{tabular}{@{}lclrrrl@{}}
\toprule
Method & Initial Guess & $\by^*$ & RRE & SSIM & PSF err. & Stop \\
\midrule
No log-det & 1 & $(0.044,0.185)$ & $0.170$ & $0.636$ & $9.552$ & Grad. \\
No log-det & 2 & $(0.068,0.311)$ & $0.170$ & $0.636$ & $9.552$ & Step \\
No log-det & 3 & $(2.560,0.059)$ & $0.163$ & $0.634$ & $3.556$ & Step \\
\addlinespace[3pt]
Marg & 1 & $(2.890,4.945)$ & $\mathbf{0.131}$ & $\mathbf{0.675}$ & $0.037$ & Step \\
Marg & 2 & $(2.890,4.945)$ & $\mathbf{0.131}$ & $\mathbf{0.675}$ & $0.037$ & Grad. \\
Marg & 3 & $(2.890,4.945)$ & $\mathbf{0.131}$ & $\mathbf{0.675}$ & $0.037$ & Grad. \\
\addlinespace[3pt]
Proposed & 1 & $(3.119,4.967)$ & $\mathbf{0.131}$ & $0.673$ & $\mathbf{0.032}$ & Step \\
Proposed & 2 & $(3.119,4.967)$ & $\mathbf{0.131}$ & $0.673$ & $\mathbf{0.032}$ & Step \\
Proposed & 3 & $(3.119,4.967)$ & $\mathbf{0.131}$ & $0.673$ & $\mathbf{0.032}$ & Step \\
\bottomrule
\end{tabular}
\end{minipage}
\end{figure}

\begin{figure}[!htbp]
\centering
\includegraphics[width=\linewidth]{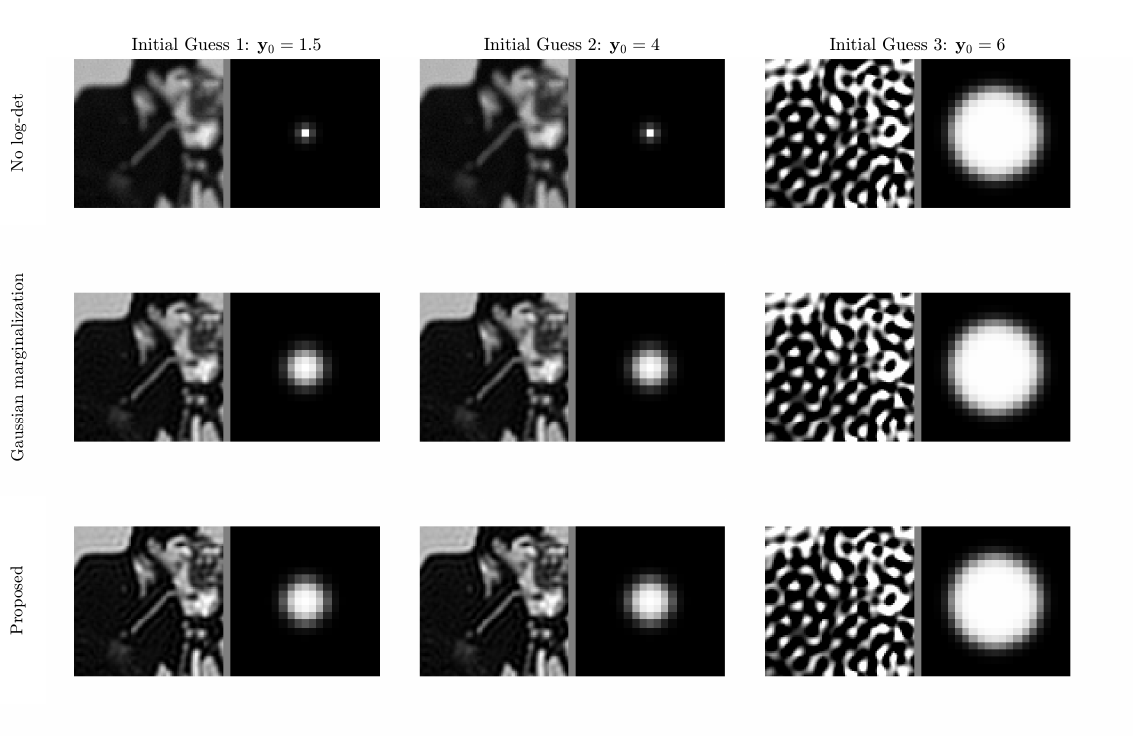}
\caption{\emph{Defocus blur: paired zoomed reconstructions and recovered PSFs from all initializations.}
Layout and display scales follow Figure~\ref{fig:composite-recon-psf-gaussian}. Columns correspond to the defocus initial guesses in Table~\ref{tab:reconstruction-defocus}.}
\label{fig:composite-recon-psf-defocus}
\par\medskip
\begin{minipage}{\linewidth}
\centering
\normalsize
\resulttablecaption{Recovery for defocus blur from all three initializations (proposed $\alpha=0.75$).}
\label{tab:reconstruction-defocus}
\setlength{\tabcolsep}{4.5pt}
\renewcommand{\arraystretch}{1.05}
\begin{tabular}{@{}lclrrrl@{}}
\toprule
Method & Initial Guess & $\by^*$ & RRE & SSIM & PSF err. & Stop \\
\midrule
No log-det & 1 & $3.68\times10^{-14}$ & $0.130$ & $0.691$ & $1.937$ & Step \\
No log-det & 2 & $4.68\times10^{-14}$ & $0.130$ & $0.691$ & $1.937$ & Step \\
No log-det & 3 & $6.057$ & $0.963$ & $0.004$ & $0.814$ & Step \\
\addlinespace[3pt]
Marg & 1 & $2.258$ & $0.104$ & $\mathbf{0.769}$ & $0.124$ & Step \\
Marg & 2 & $2.258$ & $0.104$ & $\mathbf{0.769}$ & $0.124$ & Step \\
Marg & 3 & $6.060$ & $0.964$ & $0.004$ & $0.815$ & Step \\
\addlinespace[3pt]
Proposed & 1 & $2.665$ & $\mathbf{0.103}$ & $0.753$ & $\mathbf{0.078}$ & Step \\
Proposed & 2 & $2.665$ & $\mathbf{0.103}$ & $0.753$ & $\mathbf{0.078}$ & Step \\
Proposed & 3 & $6.065$ & $0.965$ & $0.004$ & $0.815$ & Step \\
\bottomrule
\end{tabular}
\end{minipage}
\end{figure}

\begin{figure}[!htbp]
\centering
\includegraphics[width=\linewidth]{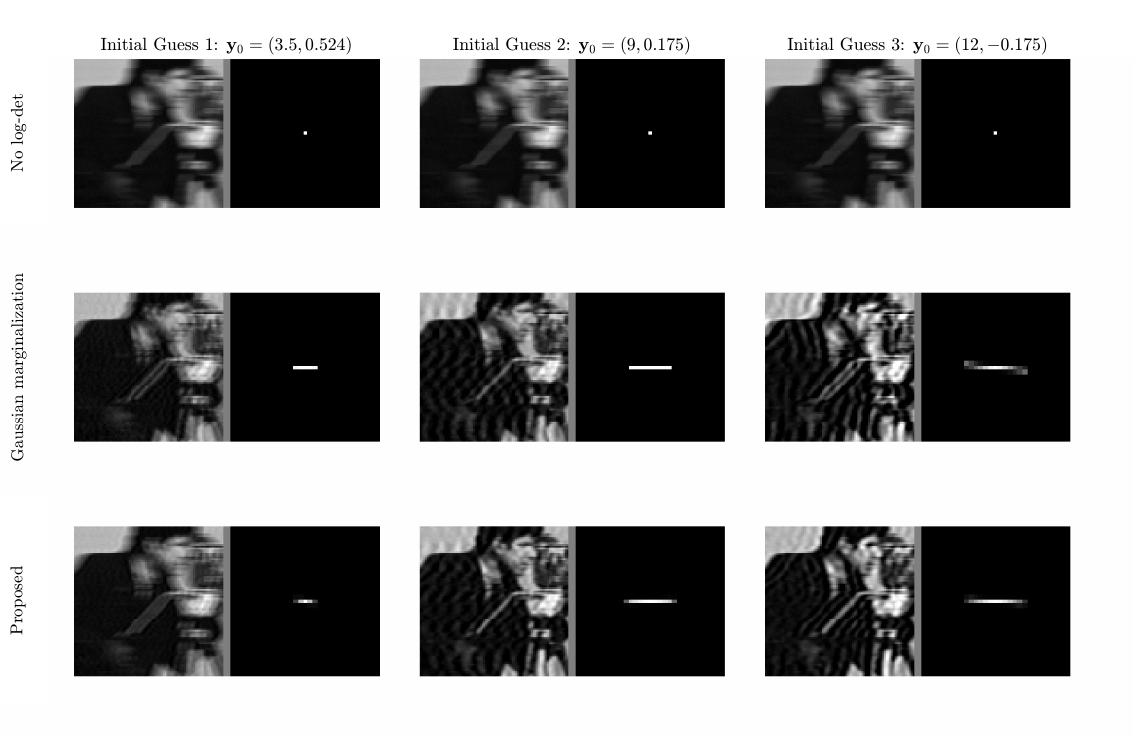}
\caption{\emph{Motion blur: paired zoomed reconstructions and recovered PSFs from all initializations.}
Layout and display scales follow Figure~\ref{fig:composite-recon-psf-gaussian}. Columns correspond to the motion initial guesses in Table~\ref{tab:reconstruction-motion}.}
\label{fig:composite-recon-psf-motion}
\par\medskip
\begin{minipage}{\linewidth}
\centering
\normalsize
\resulttablecaption{Recovery for motion blur from all three initializations (proposed $\alpha=\eta/N\approx0.5722$).}
\label{tab:reconstruction-motion}
\setlength{\tabcolsep}{4.5pt}
\renewcommand{\arraystretch}{1.05}
\begin{tabular}{@{}lclrrrl@{}}
\toprule
Method & Initial Guess & $\by^*$ & RRE & SSIM & PSF err. & Stop \\
\midrule
No log-det & 1 & $(2.34\times10^{-7},1.111)$ & $0.157$ & $\mathbf{0.723}$ & $3.092$ & Step \\
No log-det & 2 & $(4.14\times10^{-5},-3.024)$ & $0.157$ & $\mathbf{0.723}$ & $3.078$ & Step \\
No log-det & 3 & $(5.03\times10^{-5},3.077)$ & $0.157$ & $\mathbf{0.723}$ & $3.073$ & Step \\
\addlinespace[3pt]
Marg & 1 & $(4.486,0.007)$ & $0.181$ & $0.602$ & $1.000$ & Limit \\
Marg & 2 & $(9.036,-1.44\times10^{-4})$ & $0.133$ & $0.591$ & $0.253$ & Limit \\
Marg & 3 & $(13.158,-0.100)$ & $0.207$ & $0.476$ & $0.419$ & Step \\
\addlinespace[3pt]
Proposed & 1 & $(4.118,-0.232)$ & $0.158$ & $0.699$ & $1.172$ & Step \\
Proposed & 2 & $(10.018,0.052)$ & $\mathbf{0.110}$ & $0.669$ & $\mathbf{0.093}$ & Step \\
Proposed & 3 & $(12.009,-0.072)$ & $0.122$ & $0.673$ & $0.212$ & Step \\
\bottomrule
\end{tabular}
\end{minipage}
\end{figure}

We use $\beta^2=0$ in the first case and $(\beta^2,\gamma)=(1/\eta,\lambda^2)$ in the second.
For the third case, we set $\gamma=1/\eta$ and $\beta^2=\alpha/\eta$.
For each blur family, we solve the problem for
$$ \alpha\in\left\{ 0.05,\;0.1,\;0.2,\;0.35,\;0.5,\;\frac{\eta}{N}, \;0.75,\;1,\;1.5,\;2 \right\}, $$
from the same three initial guesses. We then select the value of $\alpha$ that minimizes the average image relative reconstruction error with respect to the true image across all three runs, including unsuccessful reconstructions. The selected values are $\alpha=0.5$ for Gaussian blur, $\alpha=0.75$ for defocus blur, and $\alpha=\eta/N\approx0.5722$ for motion blur. The same selected value is used for all three initial guesses within each blur family. We use the initial guesses $(4,8)$, $(2,2)$, and $(8,4)$ for Gaussian blur; $1.5$, $4$, and $6$ for defocus blur; and $(3.5,\pi/6)$, $(9,\pi/18)$, and $(12,-\pi/18)$ for motion blur. These are numbered 1-3 in Tables~\ref{tab:reconstruction-gaussian}-\ref{tab:reconstruction-motion}. We apply Algorithm~\ref{alg:HybridGN} directly to $\bff_{\beta,\gamma}(\by)$, with initial step length
$t_{\mathrm{init}}=0.25$, Armijo constant $c_1=10^{-4}$, and backtracking factor $\zeta=0.5$. We allow at most 500 iterations and use step tolerance $\epsilon_s=10^{-12}$, damping parameter $\tau_{\mathrm{base}}=10^{-6}/N$, and gradient tolerance
$\epsilon_g=10^{-3}/N$. Once the blur estimate $\by^*$ is computed, we recover the image $\bx^*=\bx(\by^*)$ using~\eqref{eq:Tikhonov-solution}.

To assess the reconstruction quality, we compute the absolute parameter error, relative residual error, image relative reconstruction error (RRE), peak signal-to-noise ratio (PSNR), and relative PSF error,
defined by
$$
\begin{aligned}
\text{Parameter error} &=\|\by^*-\by_{\mathrm{true}}\|_2,\qquad \mathrm{RResE}(\bx^*,\by^*) =\frac{ \|\bA(\by^*)\bx^*-\bb\|_2} {\|\bb\|_2}, \\
\mathrm{RRE}(\bx^*) &=\frac{\|\bx^*-\bx_{\mathrm{true}}\|_2} {\|\bx_{\mathrm{true}}\|_2},\qquad \mathrm{PSNR}(\bx^*) =10\log_{10}\!\left( \frac{x_{\max}^2} {\frac{1}{N}\|\bx^*-\bx_{\mathrm{true}}\|_2^2} \right), \\
\text{Relative PSF error} &=\frac{\|p^*-p_{\mathrm{true}}\|_2} {\|p_{\mathrm{true}}\|_2}.
\end{aligned}
$$
Here, $x_{\max}=1$, and $p_{\mathrm{true}}$ and $p^*$ denote the true and estimated PSFs, respectively.
We also report the structural similarity index (SSIM)~\cite{wang2004image}. All image quality measures are computed using the full reconstruction without clipping its intensities.

\begin{figure}[!htbp]
\centering
\includegraphics[width=\linewidth]{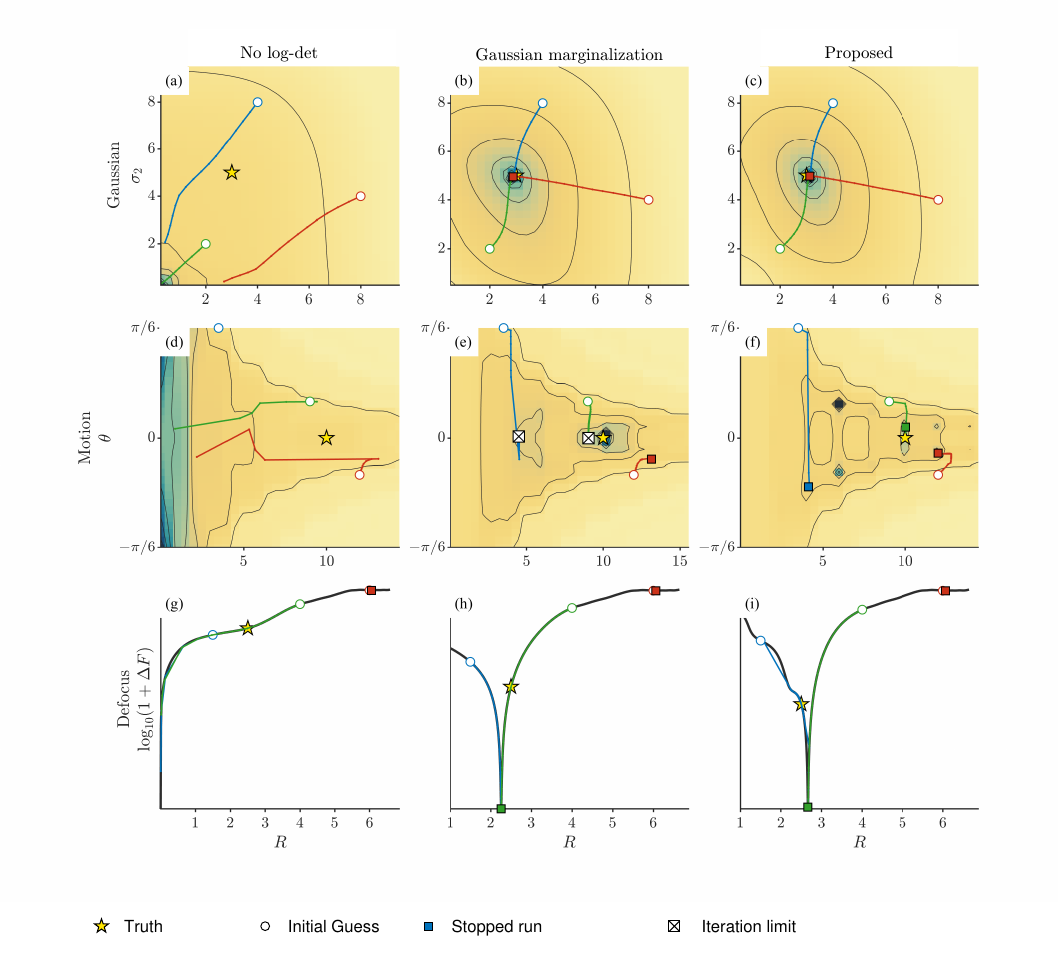}
\caption{\emph{Objective landscapes for the semi-blind deconvolution reconstruction experiments.} Rows correspond to the Gaussian, motion, and defocus PSF models, and columns are labeled ``No log-det'', ``Gaussian marginalization'', and ``Proposed''. For display, $F=N\bff_{\beta,\gamma}$ and $\Delta F$ is its value shifted by the minimum on the plotting grid; the plotted value is $\log_{10}(1+\Delta F)$. Stars mark the true parameter, open circles mark the initial guesses, curves mark optimization trajectories, filled squares mark stopped runs, and open squares with an $\times$ mark runs that reached the iteration limit.}
\label{fig:landscape-results}
\end{figure}

\medskip\noindent\textbf{Results.}
Tables~\ref{tab:reconstruction-gaussian}-\ref{tab:reconstruction-motion} report all initial guesses, including failed reconstructions; Supplementary Section~\ref{app:reconstruction-metrics} gives additional metrics. ``No log-det'' denotes $\beta^2=0$, ``Marg'' denotes Gaussian marginalization, and ``Proposed'' denotes the proposed regularization with the weights specified above. Lower image RRE and PSF error and higher SSIM indicate better reconstructions. Boldface marks each metric's best displayed value within a blur family, including ties at the displayed precision. Grad., Step, and Limit denote gradient tolerance, small-step termination, and the iteration limit. Step includes zero-step returns after failed backtracking; stopping or a small residual does not certify an accurate estimate. Figure~\ref{fig:landscape-results} shows the optimization trajectories on the corresponding landscapes.

For Gaussian blur (see Figures~\ref{fig:composite-recon-psf-gaussian} and \ref{fig:landscape-results}), without the log-determinant regularizer, both estimated widths degenerate from the first two initial guesses, and one degenerates from the third. Both proposed log-det regularizers recover nondegenerate estimates from all three initial guesses and yield similar parameter and image reconstruction errors, whereas the auxiliary choice \(\tilde\bL=\bI\) yields a smaller relative PSF error (see Tables~\ref{tab:reconstruction-gaussian} and \ref{tab:additional-gaussian}).

For defocus blur (see Figures~\ref{fig:composite-recon-psf-defocus} and \ref{fig:landscape-results}), both log-det regularizers improve estimates from the first two initial guesses, with the proposed regularizer ($\tilde\bL=\bI$) giving smaller parameter, image RRE, and PSF errors than marginalization. All three methods remain in the large-radius region from the third initial guess, with image RRE near $0.96$ (see Tables~\ref{tab:reconstruction-defocus} and \ref{tab:additional-defocus}).

For motion blur (Figures~\ref{fig:composite-recon-psf-motion} and \ref{fig:landscape-results}), the proposed regularizer improves parameter, image, and PSF estimates over marginalization from the second and third initial guesses (see Tables~\ref{tab:reconstruction-motion} and \ref{tab:additional-motion}). From the first, it lowers image RRE but raises parameter and PSF errors. Without the log-determinant regularizer, the estimated motion lengths are near the origin from all three initial guesses, despite achieving the highest SSIM. Marginalization reaches the maximum number of iterations from its first two initial guesses and stops at an inaccurate estimate from the third.

These comparisons use determinant weights selected with the true image;
Section~\ref{ssec:practical-parameter-choice} examines selection from
the observations.

\subsection{Sensitivity analysis} \label{ssec:auxiliary-prior-sensitivity}
\begin{figure}[!htbp]
\centering

\includegraphics[width=\linewidth]{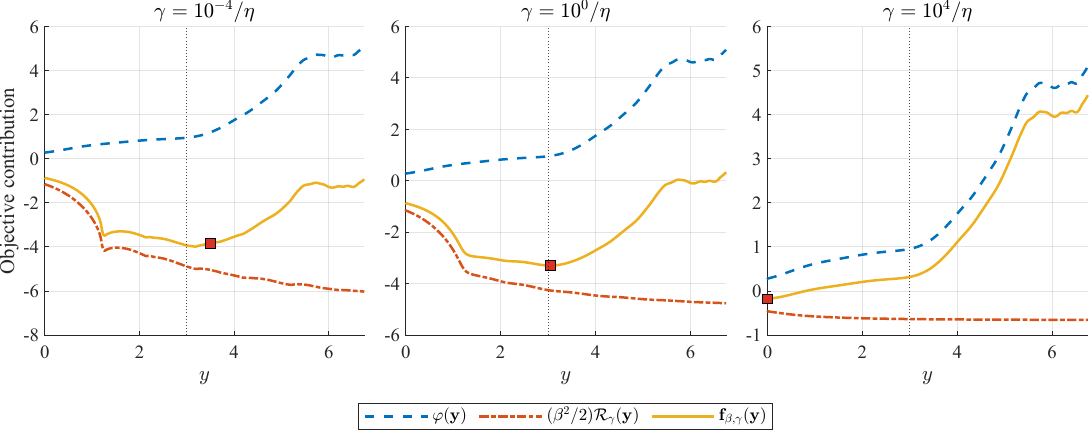}

\medskip

\includegraphics[width=\linewidth]{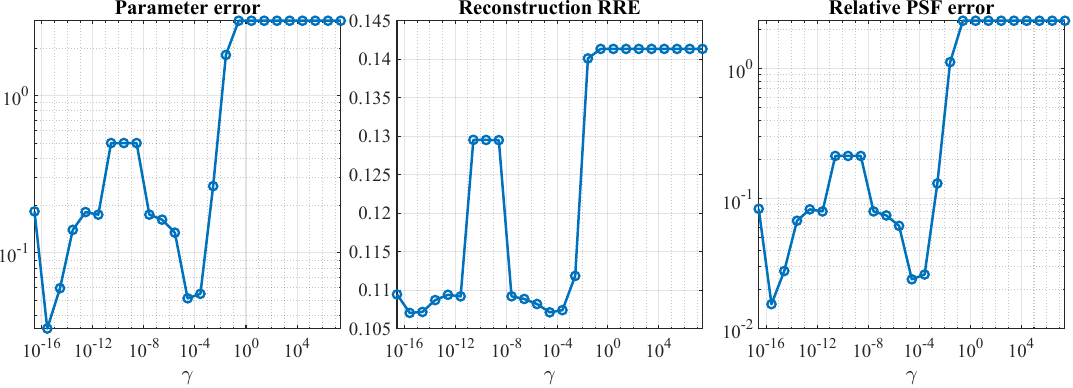}

\caption{\emph{Auxiliary-scale sensitivity for smooth defocus blur, with $\beta^2=1/N$ and $\tilde\bL=\bI$.} Top: the reduced Tikhonov term $\varphi(\by)$, the weighted determinant penalty $\frac{\beta^2}{2}\mathcal R_\gamma(\by)$, and the total objective $\bff_{\beta,\gamma}(\by)$. From left to right, the panels use $\gamma=10^{-4}/\eta$, $1/\eta$, and $10^4/\eta$. The vertical dotted line marks the true radius, and the square marks the computed estimate. Bottom: parameter error, image RRE, and relative PSF error across all 25 scales.}
\label{fig:delta-sensitivity}
\end{figure}

In this experiment, we vary the auxiliary scale $\gamma$ with $\tilde\bA=\bA$ and $\tilde\bL=\bI$ fixed, isolating the determinant's spectral floor from image regularization and examining the behavior in Proposition~\ref{prop:delta_limits}.
We use the same \emph{Cameraman} image, smooth defocus model, $1\%$ relative Gaussian noise, and periodic second-difference image regularizer, with $\lambda^2=100/\eta$. Here the true radius is $\by_{\mathrm{true}}=3$ and the initial guess is $\by_0=3.5$. The solver settings are unchanged except for the step tolerance, which is $10^{-10}$.
We minimize $\bff_{\beta,\gamma}(\by)$ in~\eqref{eq:implementation-objective} with the fixed weight $\beta^2=1/N$, where $N=256^2$. Only $\gamma$ varies over the 25 scales
$$ \gamma=10^j/\eta,\qquad j=-12,-11,\ldots,12. $$
For each scale, we plot the parameter error, image RRE, and relative PSF error in Figure~\ref{fig:delta-sensitivity}. The PSF error compares the estimated and true smooth defocus PSFs~\eqref{eq:defocus-psf}. Figure~\ref{fig:delta-sensitivity} also shows the contributions of $\varphi(\by)$ and the regularization terms in $\bff_{\beta,\gamma}(\by)$, separately, for $\gamma=10^{-4}/\eta,\allowbreak 1/\eta,\allowbreak 10^4/\eta$.

\medskip\noindent\textbf{Results.}
At small $\gamma$, estimates vary irregularly, and the determinant responds strongly to small eigenvalues of $\bA(\by)^\top\bA(\by)$, consistent with the sharp slope changes at $\gamma=10^{-4}/\eta$. At $\gamma=1/\eta$, the estimated radius is $3.051$, with image RRE equal to $0.107$ and relative PSF error $0.0239$. This intermediate scale gives accurate recovery in the present example, although some smaller scales yield lower errors. The sweep shows that increasing the auxiliary scale changes both the objective's variation with the blur parameter and the recovered estimate; its effect on reconstruction accuracy is not monotone. For large $\gamma$, the auxiliary term dominates $\tilde\bG_\gamma(\by)$ and the determinant becomes nearly constant in $\by$, as predicted by~\eqref{eq:delta-large-regime}. Up to terms independent of $\by$, this limiting objective retains image regularization but lacks regularization of the blur parameter. When $\bA(\by)\to\bI$, the no-blur limit is characterized by~\eqref{eq:degnerate-objective}, with the corresponding reconstruction given in~\eqref{eq: degenerate-tikhonov-solution}. In the present smoothed-defocus experiment, the estimated radii approach zero at scales of $10^4/\eta$ and larger, with image RRE approximately $0.141$. This behavior is consistent with a preference for low-blur solutions, although the fixed smoothing prevents the defocus PSF from reaching the exact identity limit.

\subsection{An illustration of automatic parameter selection}
\label{ssec:practical-parameter-choice}

The preceding experiments prescribe regularization parameters to study their effects on the estimates. We now give a preliminary, single-sample illustration of choosing these parameters from the observations: we estimate $\eta$ using the median absolute deviation (MAD) of image differences, select $\lambda$ using a conditional form of generalized cross-validation (GCV), and select the prior exponent $\alpha$ through a bounded search. This use of GCV is motivated by nonlinear local-linearization methods~\cite{osullivan1985cross,haber2000gcv} and regularized VarPro methods~\cite{chung2010efficient, Espanol_2023}. The score holds the blur fixed, omitting its fitted sensitivity to the data.

We use again the $256\times256$ \emph{Cameraman} image with Gaussian PSF with parameters $\by_{\rm true}=(3,2)$, $5\%$ relative noise, and initial guess $\by_0=(4,4)$. The observed data are shown in the last panel of Figure~\ref{fig:blurred-data}. We use periodic first-difference image regularization. We fix $\tilde\bA=\bA$, $\tilde\bL=\bI$, and $\gamma=1/\eta_{\rm MAD}$, with $\beta^2=\alpha/\eta_{\rm MAD}$ in~\eqref{eq:implementation-objective}. The selected fit and $\alpha=0$ control share the observation, noise estimate, and conditional-GCV selection of $\lambda$; every trial starts from an initial guess $\by_0$. Supplementary Section~\ref{ssec:profiled-gcv-beta} gives the procedure, search bounds, and stopping tolerances.

\begin{table}[!htbp]
\centering\small
\caption{Single-sample automatic parameter selection for Gaussian blur. Ground truth is used only to compute the reported recovery errors.}
\label{tab:auto}
\begin{tabular}{lrrrrrr}
\toprule
Method & $\alpha$ & $\lambda$ & $\by^*$ & RRE & PSF error & SSIM\\
\midrule
No log-det & 0.000 & $0.409$ & $(0.160,0.374)$ & 0.139 & 5.881 & 0.621 \\
Automatic selection & 0.331 & $0.123$ & $(3.437,2.035)$ & 0.119 & 0.119 & 0.589 \\
\bottomrule
\end{tabular}
\end{table}

\begin{figure}[!htbp]
\centering
\includegraphics[width=\linewidth]{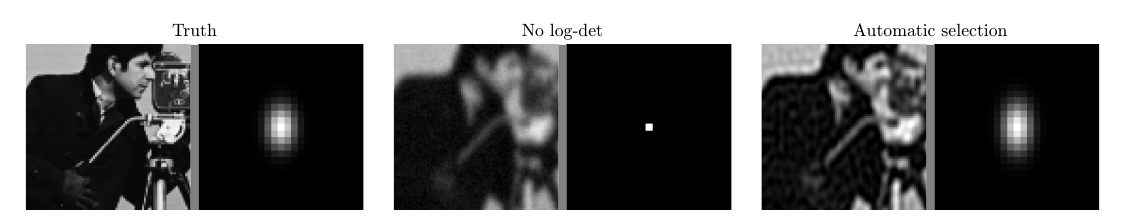}
\caption{\emph{Single-sample automatic parameter selection.} Each pair shows a zoomed reconstruction (left) and its PSF (right). Pairs show the truth, the no log-det control, and the automatically selected fit. Reconstructions use the same image crop as Figures~\ref{fig:composite-recon-psf-gaussian}-\ref{fig:composite-recon-psf-motion} and the common intensity scale $[0,1]$. Central PSF crops use independent color scales; full kernels sum to one.}
\label{fig:gcv-single-sample}
\end{figure}

\medskip\noindent\textbf{Results.}
Table~\ref{tab:auto} and Figure~\ref{fig:gcv-single-sample} show that automatic selection avoids the control's near-degenerate blur estimate and improves image RRE and PSF recovery. This comparison assesses the full selection procedure, since the
image regularization parameter $\lambda$ is reselected alongside
$\alpha$. The sharper reconstruction has lower SSIM and amplified
texture. The estimated noise precision is $0.981$ times the reference precision, and the bounded search and both final alternating fits meet their stopping criteria. Further work is needed to evaluate the procedure across images and noise realizations.

\section{Conclusions}
\label{sec: Conclusions}

Motivated by Gaussian marginalization, we have proposed a family of log-determinant regularizers for the nonlinear variables in separable nonlinear least squares problems and solved the resulting RSNLS problems using a variable projection method. We have examined how these regularizers affect the objective landscape and the behavior of the numerical method, and evaluated their performance in semi-blind image deblurring.

Our analysis clarifies the roles of the auxiliary operators and the determinant weight. The auxiliary operators determine how the regularizer varies with the blur parameters, while the weight controls its contribution relative to the reduced Tikhonov term. Under the stated spectral assumptions, a sufficiently large weight excludes the no-blur limit as a local minimizer. The auxiliary-scale analysis also shows that the leading terms in the large-scale expansion are independent of the blur parameters. These findings highlight the importance of choosing both the auxiliary operators and the regularization weight. The numerical experiments show that these choices affect both the objective landscape and the reconstructed image and blur. With weights selected using the true image, the proposed regularizer improves recovery over Gaussian marginalization in several settings. However, some runs remain sensitive to the initial guess and converge to inaccurate reconstructions. Further work should examine how the auxiliary choices affect both the minimizers of the reduced objective and the convergence of the numerical method.

Selecting the regularization parameters using only the observations remains an important practical question. The single-sample experiment with conditional generalized cross-validation provides an initial illustration of automatic selection. Further work should investigate the joint selection of the determinant weight, the auxiliary scale, and the image regularization parameter, including criteria that account for the dependence of the estimated blur on the data. Experiments with different images, noise realizations, initial guesses, and discrepancies between the assumed and true blur models would help assess the reliability of the approach. Extensions beyond periodic boundary conditions and quadratic image regularization could also be explored using stochastic log-determinant approximations~\cite{ubaru2017fast} and nonsmooth variable projection~\cite{van2021variable}.

\section*{Acknowledgments}

The authors acknowledge using Claude, ChatGPT, and Codex (OpenAI) to assist with manuscript organization and revision, \LaTeX{} preparation, and the implementation and debugging of numerical code.
These tools were also consulted during the general research and development of this manuscript. The authors assume full responsibility for the accuracy and integrity of all content in this paper.

J. T. Dworaczyk is supported by the Department of War’s Science, Mathematics, and Research for Transformation (SMART) Scholarship-for-Service Program, funded by: OUSD/R\&E (The Under Secretary of Defense-Research and Engineering), National Defense Education Program (NDEP) / BA-1, Basic Research. Any opinions, findings, and conclusions or recommendations expressed in this material are those of the authors and do not necessarily reflect the views of the Department of War.

\bibliographystyle{plain}
\bibliography{biblio}

@article{chen2018regularized,
  title={A regularized variable projection algorithm for separable nonlinear least-squares problems},
  author={Chen, G.-Y. and Gan, M. and Chen, C. L. P. and Li, H.-X.},
  journal={IEEE Transactions on Automatic Control},
  volume={64},
  number={2},
  pages={526--537},
  year={2019},
  publisher={IEEE}
}

@article{fu2019tikhonov,
  title={Tikhonov regularized variable projection algorithms for separable nonlinear least squares problems},
  author={Fu, Z. and Guo, L.},
  journal={Complexity},
  volume={2019},
  number={1},
  pages={4861708},
  year={2019},
  publisher={Hindawi}
}

@article{golub1973differentiation,
  title={The differentiation of pseudo-inverses and nonlinear least squares problems whose variables separate},
  author={Golub, G. H. and Pereyra, V.},
  journal={SIAM Journal on Numerical Analysis},
  volume={10},
  number={2},
  pages={413--432},
  year={1973},
  publisher={SIAM}
}

@article{cornelio2014constrained,
  title={Constrained numerical optimization methods for blind deconvolution},
  author={Cornelio, A. and Loli Piccolomini, E. and Nagy, J. G.},
  journal={Numerical Algorithms},
  volume={65},
  number={1},
  pages={23--42},
  year={2014},
  doi={10.1007/s11075-013-9693-z},
  publisher={Springer}
}

@article{chung2010efficient,
  title={An efficient iterative approach for large-scale separable nonlinear inverse problems},
  author={J. Chung and J. G. Nagy},
  journal={SIAM Journal on Scientific Computing},
  volume={31},
  number={6},
  pages={4654--4674},
  year={2010},
  publisher={SIAM}
}

@article{Espanol_2023,
year = {2023},
publisher = {IOP Publishing},
volume = {39},
number = {8},
pages = {084002},
author = {M. I. Espa{\~n}ol and M. Pasha},
title = {Variable projection methods for separable nonlinear inverse problems with general-form {T}ikhonov regularization},
journal = {Inverse Problems}}

@article{gazzola2021regularization,
  title={Regularization by inexact {K}rylov methods with applications to blind deblurring},
  author={Gazzola, S. and Sabat{\'e} Landman, M.},
  journal={SIAM Journal on Matrix Analysis and Applications},
  volume={42},
  number={4},
  pages={1528--1552},
  year={2021},
  publisher={SIAM}
}

@article{wang2004image,
  title={Image quality assessment: from error visibility to structural similarity},
  author={Wang, Z. and Bovik, A. C. and Sheikh, H. R. and Simoncelli, E. P.},
  journal={IEEE Transactions on Image Processing},
  volume={13},
  number={4},
  pages={600--612},
  year={2004},
  publisher={IEEE}
}

@inproceedings{levin2011efficient,
  title={Efficient marginal likelihood optimization in blind deconvolution},
  author={Levin, A. and Weiss, Y. and Durand, F. and Freeman, W. T.},
  booktitle={CVPR 2011},
  pages={2657--2664},
  year={2011},
  organization={IEEE}
}

@inproceedings{perrone2014total,
  title={Total variation blind deconvolution: the devil is in the details},
  author={Perrone, D. and Favaro, P.},
  booktitle={Proceedings of the IEEE Conference on Computer Vision and Pattern Recognition},
  pages={2909--2916},
  year={2014}
}

@book{rue2005gaussian,
  title={{G}aussian {M}arkov random fields: theory and applications},
  author={Rue, H. and Held, L.},
  year={2005},
  publisher={Chapman and Hall/CRC}
}

@book{tenorio2017data,
  title={An Introduction to Data Analysis and Uncertainty Quantification for Inverse Problems},
  author={Tenorio, L.},
  year={2017},
  publisher={Society for Industrial and Applied Mathematics},
  address={Philadelphia, PA},
  doi={10.1137/1.9781611974928}
}

@article{van2021variable,
  title={Variable projection for nonsmooth problems},
  author={T. Van Leeuwen and A. Y. Aravkin},
  journal={SIAM Journal on Scientific Computing},
  volume={43},
  number={5},
  pages={S249--S268},
  year={2021},
  publisher={SIAM}
}

@book{nocedal2006numerical,
  title={Numerical Optimization},
  author={Nocedal, J. and Wright, S. J.},
  year={2006},
  edition={second},
  publisher={Springer Science \& Business Media}
}

@article{pereyra2019imaging,
  title={Imaging applications with variable projections},
  author={Pereyra, V. and Scherer, G.},
  journal={American Journal of Computational Mathematics},
  volume={9},
  number={04},
  pages={261--281},
  year={2019}
}

@book{deblurring,
 author = {P. C. Hansen and J. G. Nagy and D. P. O'Leary},
 title = {Deblurring {I}mages: {M}atrices, {Spectra} and {Filtering}},
 year = {2006},
  publisher ={SIAM, Philadelphia},
 }

@article{espanol2025local,
  title={Local Convergence Analysis of a Variable Projection Method for Regularized Separable Nonlinear Inverse Problems},
  author={M. I. Espa{\~n}ol and G. Jeronimo},
  journal={SIAM Journal on Matrix Analysis and Applications},
  volume={46},
  number={2},
  pages={858--878},
  year={2025},
  publisher={SIAM}
}

@article{marmin2025varprox,
  title={{VARPROX}: a primal-dual variable projection method for the minimization of penalized separable non-linear least squares},
  author={Marmin, A. and Richard, F. J. P.},
  journal={Inverse Problems},
  volume={41},
  number={10},
  pages={105006},
  year={2025},
  doi={10.1088/1361-6420/ae0e48},
  publisher={IOP Publishing}
}

@inproceedings{levin2009understanding,
  title={Understanding and evaluating blind deconvolution algorithms},
  author={Levin, A. and Weiss, Y. and Durand, F. and Freeman, W. T.},
  booktitle={2009 IEEE conference on computer vision and pattern recognition},
  pages={1964--1971},
  year={2009},
  organization={IEEE}
}

@article{salzer2026variable,
  title={Variable Projection Methods for Solving Regularized Separable Inverse Problems with Applications to Semi-Blind Image Deblurring},
  author={Comerso Salzer, D. B. and Espa{\~n}ol, M. I. and Jeronimo, G.},
  journal={arXiv preprint arXiv:2601.05224},
  year={2026}
}

@article{mbakam2024marginal,
  title={Marginal likelihood estimation in semiblind image deconvolution: A stochastic approximation approach},
  author={Mbakam, C. K. and Pereyra, M. and Giovannelli, J.},
  journal={SIAM Journal on Imaging Sciences},
  volume={17},
  number={2},
  pages={1206--1254},
  year={2024},
  publisher={SIAM}
}

@article{golub1979gcv,
 author={Golub, G. H. and Heath, M. and Wahba, G.},
 title={Generalized Cross-Validation as a Method for Choosing a Good Ridge Parameter},
 journal={Technometrics}, volume={21}, number={2}, pages={215--223}, year={1979},
 doi={10.1080/00401706.1979.10489751}
}

@article{osullivan1985cross,
  author={O'Sullivan, F. and Wahba, G.},
  title={A Cross Validated {Bayesian} Retrieval Algorithm for Nonlinear Remote Sensing Experiments},
  journal={Journal of Computational Physics},
  volume={59},
  number={3},
  pages={441--455},
  year={1985},
  doi={10.1016/0021-9991(85)90121-4}
}

@article{haber2000gcv,
  author={Haber, E. and Oldenburg, D.},
  title={A {GCV} Based Method for Nonlinear Ill-Posed Problems},
  journal={Computational Geosciences},
  volume={4},
  number={1},
  pages={41--63},
  year={2000},
  doi={10.1023/A:1011599530422}
}

@article{sweeney2025parameter,
  author={Sweeney, B. and Renaut, R. and Espa{\~n}ol, M. I.},
  title={Parameter Selection by {GCV} and a {$\chi^2$} Test within Iterative Methods for {$\ell_1$}-Regularized Inverse Problems},
  journal={SIAM Journal on Scientific Computing},
  pages={S112--S134},
  year={2025},
  doi={10.1137/24M165805X}
}

@article{ubaru2017fast,
 author={Ubaru, S. and Chen, J. and Saad, Y.},
 title={Fast Estimation of {$\operatorname{tr}(f(A))$} via Stochastic {Lanczos} Quadrature},
 journal={SIAM Journal on Matrix Analysis and Applications},
 volume={38}, number={4}, pages={1075--1099}, year={2017},
 doi={10.1137/16M1104974}
}

@article{fergus2006removing,
  author  = {Fergus, R. and Singh, B. and Hertzmann, A.
             and Roweis, Sam T. and Freeman, William T.},
  title   = {Removing Camera Shake from a Single Photograph},
  journal = {ACM Transactions on Graphics},
  volume  = {25},
  number  = {3},
  pages   = {787--794},
  year    = {2006},
  doi     = {10.1145/1141911.1141956}
}

@article{wipf2014revisiting,
  author  = {Wipf, D. and Zhang, H.},
  title   = {Revisiting {Bayesian} Blind Deconvolution},
  journal = {Journal of Machine Learning Research},
  volume  = {15},
  number  = {111},
  pages   = {3775--3814},
  year    = {2014},
  url     = {https://www.jmlr.org/papers/v15/wipf14a.html}
}

@article{buser2026majorization,
  title={A Majorization-Minimization with Monte Carlo Approach for Hyperparameter Estimation},
  author={Buser, E. and Chung, J. and D{\'\i}az, H. and Saibaba, A. K.},
  journal={arXiv preprint arXiv:2605.13620},
  year={2026}
}

\renewcommand{\appendixname}{}

\setcounter{section}{0}
\renewcommand{\thesection}{SM\arabic{section}}

\setcounter{equation}{0}
\setcounter{figure}{0}
\setcounter{table}{0}
\setcounter{algorithm}{0}
\counterwithin{equation}{section}
\renewcommand{\thefigure}{SM\arabic{figure}}
\renewcommand{\thetable}{SM\arabic{table}}
\renewcommand{\thealgorithm}{SM\arabic{algorithm}}
\renewcommand{\theHsection}{SM\arabic{section}}
\renewcommand{\theHequation}{\theHsection.\arabic{equation}}
\renewcommand{\theHfigure}{SM\arabic{figure}}
\renewcommand{\theHtable}{SM\arabic{table}}
\renewcommand{\theHalgorithm}{SM\arabic{algorithm}}

\bigskip
\begin{center}
\large\bfseries SUPPLEMENTARY MATERIALS
\end{center}
\medskip

Sections, equations, tables, and algorithms with the prefix SM belong
to this supplementary material. Unprefixed references to these items
refer to the main manuscript.

\section{Derivative formulas for the log-determinant term}
\label{app:logdet-derivatives}
We derive the matrix identities underlying the Fourier formulas in
Section~\ref{ssec:fft-evaluation}, using the same derivative notation.
Consider the auxiliary matrix
$\tilde\bG_\gamma$ in~\eqref{eq:auxiliary-gram}, with fixed $\gamma$
and $\tilde\bL$, and assume that $\tilde\bA$ is twice
continuously differentiable and $\tilde\bG_\gamma$ is positive definite.
For brevity, write
$\tilde\bG_{\gamma,i}=\partial_{y_i}\tilde\bG_\gamma$ and
$\tilde\bG_{\gamma,i\ell}=\partial^2_{y_i y_\ell}\tilde\bG_\gamma$.
The product rule gives
$$
\begin{aligned}
\tilde\bG_{\gamma,i} &=\left[ (\partial_{y_i}\tilde\bA)^\top\tilde\bA +\tilde\bA^\top\partial_{y_i}\tilde\bA\right], \\
\tilde\bG_{\gamma,i\ell} &=\left[ (\partial^2_{y_i y_\ell}\tilde\bA)^\top\tilde\bA +(\partial_{y_i}\tilde\bA)^\top\partial_{y_\ell}\tilde\bA +(\partial_{y_\ell}\tilde\bA)^\top\partial_{y_i}\tilde\bA \right. \\
&\qquad\left. +\tilde\bA^\top\partial^2_{y_i y_\ell}\tilde\bA\right].
\end{aligned}
$$
Jacobi's formula yields
$$ \partial_{y_i}\left(\frac12\log\det\tilde\bG_\gamma\right) =\frac12\operatorname{tr}(\tilde\bG_\gamma^{-1}\tilde\bG_{\gamma,i}) =\operatorname{tr}\left( \tilde\bG_\gamma^{-1}\tilde\bA^\top\partial_{y_i}\tilde\bA\right), $$
where the last equality follows by transposing one product term
and using cyclic invariance of the trace.
Differentiating and using
$\partial_{y_\ell}\tilde\bG_\gamma^{-1}
=-\tilde\bG_\gamma^{-1}\tilde\bG_{\gamma,\ell}\tilde\bG_\gamma^{-1}$ gives
$$ \partial^2_{y_i y_\ell}\left(\frac12\log\det\tilde\bG_\gamma\right) =\frac12\operatorname{tr}\left( \tilde\bG_\gamma^{-1}\tilde\bG_{\gamma,i\ell} -\tilde\bG_\gamma^{-1}\tilde\bG_{\gamma,\ell}\tilde\bG_\gamma^{-1}\tilde\bG_{\gamma,i} \right). $$
For the periodic operators used in Section~\ref{ssec:fft-evaluation},
substitution of their Fourier diagonalization reduces these traces
to the componentwise formulas given there.

\section{Gauss--Newton approximation of the VarPro term}
\label{app:varpro-gn}

Using the \linebreak Fourier conventions of Section~\ref{sec: analysis}, define the Fourier residual vector $\boldsymbol\rho(\by)$ by
\begin{equation}\label{eq:compressed-residual}
\rho_k(\by):= \frac{\lambda|\widehat s_k||\widehat b_k|} {\sqrt{|\widehat c_k(\by)|^2+\lambda^2|\widehat s_k|^2}}.
\end{equation}
Then $\tfrac12\|\boldsymbol\rho(\by)\|_2^2=\varphi(\by)$.
For fixed $\lambda$, differentiation gives
\begin{equation}\label{eq:compressed-jacobian}
(\bJ_\rho)_{ki}=- \frac{\lambda|\widehat s_k||\widehat b_k| \operatorname{Re}(\overline{\widehat c_k}\, \partial_{y_i}\widehat c_k)} {(|\widehat c_k|^2+\lambda^2|\widehat s_k|^2)^{3/2}}.
\end{equation}
Consequently, the exact gradient is $\bJ_\rho^\top\boldsymbol\rho$,
and the Gauss--Newton matrix is $\bJ_\rho^\top \bJ_\rho$, giving
\eqref{eq:fft-reduced-gradient} and~\eqref{eq:fft-gn-matrix}.
The exact Hessian is
\begin{equation}\label{eq:exact-residual-hessian}
\nabla^2\left(\frac12\|\boldsymbol\rho\|_2^2\right) =\bJ_\rho^\top \bJ_\rho+\sum_k\rho_k\nabla^2\rho_k.
\end{equation}
The Gauss--Newton approximation omits the second term
in~\eqref{eq:exact-residual-hessian}.

\section{Data-driven parameter selection}
\label{ssec:profiled-gcv-beta}

We specify the procedure used in Section~\ref{ssec:practical-parameter-choice} and Algorithm~\ref{alg:profiled-gcv-beta}; numerical settings follow the algorithm.

\medskip\noindent\emph{Motivation and precedent.}
Precedents include O'Sullivan and Wahba's~\cite{osullivan1985cross} approximate criterion based on a linearized influence matrix after a nonlinear fit at each candidate regularization parameter, and Haber and Oldenburg's~\cite{haber2000gcv} GCV updates for linearized problems within nonlinear iterations. Related updates occur in regularized VarPro~\cite{chung2010efficient,Espanol_2023}. Sweeney, Renaut, and Espa\~nol~\cite{sweeney2025parameter} derive inner Tikhonov GCV criteria for iterative $\ell_1$ regularization with a fixed forward operator. These precedents motivate our conditional score but do not establish a predictive-risk interpretation for the complete estimator.

\medskip\noindent\emph{Conditional GCV selection of $\lambda$.}
For fixed $\by$ and $\lambda$, the Tikhonov estimate $\bx^*(\by,\lambda)$ in~\eqref{eq:Tikhonov-solution} has filter factors and effective degrees of freedom
\begin{equation}\label{eq:gcv-filter-factors}
h_k(\by,\lambda) = \frac{|\widehat c_k(\by)|^2} {|\widehat c_k(\by)|^2+\lambda^2|\widehat s_k|^2}, \qquad \operatorname{df}(\by,\lambda)=\sum_{k=0}^{n-1} h_k(\by,\lambda).
\end{equation}
For fixed $\by$, the degrees of freedom in~\eqref{eq:gcv-filter-factors} give the usual linear GCV functional~\cite{golub1979gcv}
\begin{equation}\label{eq:conditional-gcv}
\operatorname{GCV}(\by,\lambda) = \frac{n^{-1}\|\bb-\bA(\by)\bx^*(\by,\lambda)\|_2^2} {\left(1-\operatorname{df}(\by,\lambda)/n\right)^2},
\end{equation}
and we define the bounded conditional choice by
\begin{equation}\label{eq:lambda-gcv}
\lambda_{\rm GCV}(\by) \in\argmin_{\lambda\in[10^{-8},10^3]}\operatorname{GCV}(\by,\lambda).
\end{equation}

\medskip\noindent\emph{Noise-precision estimate.}
For the observed image $\bB\in\bbR^{m\times m}$, form the diagonal second differences for $1\leq i,j<m$:
$$ d_{ij}=[\bB]_{i+1,j+1}-[\bB]_{i,j+1}-[\bB]_{i+1,j}+[\bB]_{i,j}. $$
The stencil has Euclidean norm two. The normalized median absolute deviation (MAD)
estimate~\cite[p.~29]{tenorio2017data} is
\begin{equation}\label{eq:difference-mad}
\sigma_{\rm MAD} = \frac{\operatorname{median}_{i,j} |d_{ij}-\operatorname{median}(d)|} {2\,\Phi^{-1}(3/4)}, \qquad \eta_{\rm MAD}=\sigma_{\rm MAD}^{-2},
\end{equation}
where $\Phi$ is the standard-normal distribution function.

\medskip\noindent\emph{Inner alternating problem at fixed $\alpha$.}
With $\tilde\bA=\bA$ and $\tilde\bL=\bI$, we use~\eqref{eq:implementation-objective} with
$$ \beta^2=\alpha/\eta_{\rm MAD},\qquad \gamma=1/\eta_{\rm MAD}. $$
Local solves use the deterministic objective $\bff_{\beta,\gamma}$. For fixed $\alpha$, starting from initial guess $\by^{(0)}$, we alternate
\begin{equation}\label{eq:gcv-beta-alternation}
\lambda^{(k)}=\lambda_{\rm GCV}(\by^{(k)}), \qquad \by^{(k+1)} \approx\argmin_{\by\in\calY}\bff_{\beta,\gamma}(\by;\lambda^{(k)}).
\end{equation}
Each blur update in~\eqref{eq:gcv-beta-alternation} holds $\lambda^{(k)}$ fixed and uses $\by^{(k)}$ as its initial guess; $\approx\argmin$ denotes the returned local iterate. Algorithm~\ref{alg:profiled-gcv-beta} returns the final $\by^*_\alpha$ and $\lambda^*_\alpha$.

\medskip\noindent\emph{Outer continuous selection of $\alpha$.}
At the returned inner fit, define
\begin{equation}\label{eq:profiled-gcv-beta}
\mathcal G(\alpha) :=\operatorname{GCV} (\by^*_\alpha,\lambda^*_\alpha),
\end{equation}
and select
\begin{equation}\label{eq:beta-gcv-choice}
\alpha^* \approx\argmin_{0\leq\alpha\leq1}\mathcal G(\alpha).
\end{equation}
We approximately solve~\eqref{eq:beta-gcv-choice} with \texttt{fminbnd} in \textsc{Matlab} and use its returned interior candidate. The $\alpha=0$ control is evaluated separately, not included in an endpoint comparison.

\medskip\noindent\emph{Scope of the conditional criterion.}
Using the fixed-blur linear criterion~\eqref{eq:conditional-gcv} after nonlinear fitting gives a GCV-style selection rule. If
$\bS(\by,\lambda)=\bA(\by)(\bA(\by)^\top\bA(\by)+\lambda^2\bL^\top\bL)^{-1}\bA(\by)^\top$, then, for fixed $\lambda$ and differentiable $\by^*(\bb)$,
\begin{equation}\label{eq:conditional-trace-limit}
D_\bb[\bS(\by^*(\bb),\lambda)\bb] =\bS+\sum_{i=1}^r (\partial_{y_i}\bS)\bb\,(D_\bb\,y_i^*).
\end{equation}
Here $D_\bb$ denotes the Jacobian with respect to the data vector $\bb$, rather than a derivative with respect to a blur parameter. The conditional trace omits the second, fitted-blur sensitivity term in~\eqref{eq:conditional-trace-limit}. Selecting $\lambda$ and $\alpha$ and estimating $\eta$ add further data dependence. This is a practical parameter-choice rule; a risk analysis for the complete procedure remains separate.

\begin{algorithm}
\caption{Data-driven parameter selection with difference MAD, nested GCV, and regularized VarPro}
\label{alg:profiled-gcv-beta}
\begin{algorithmic}[1]
\Require Observed image $\bB$, initial PSF parameter $\by^{(0)}$, image-regularizer symbol $\widehat s$, stopping tolerances $\epsilon_y,\epsilon_\lambda$, and alternation budget $K_{\rm alt}$
\State Estimate $\eta_{\rm MAD}$ from the diagonal differences using~\eqref{eq:difference-mad}
\State Set $\tilde\bA=\bA$ and $\tilde\bL=\bI$
\Function{NestedGCVScore}{$\alpha$}
    \State $\beta\gets\sqrt{\alpha/\eta_{\rm MAD}}$, $\gamma\gets1/\eta_{\rm MAD}$
    \State $\by\gets\by^{(0)}$, $\lambda_{\rm old}\gets\infty$
    \Repeat
        \State $\lambda\gets\texttt{fminbnd}(\operatorname{GCV}(\by,\cdot),10^{-8},10^3)$ using~\eqref{eq:conditional-gcv} and~\eqref{eq:lambda-gcv}
        \State $\by_{\rm new}\gets\text{local minimization of }\bff_{\beta,\gamma}(\cdot;\lambda)$ using Algorithm~\ref{alg:HybridGN}, with $\lambda$ fixed
        \State $r_y\gets\|\by_{\rm new}-\by\|_2/\max\{1,\|\by\|_2\}$
        \State $r_\lambda\gets|\lambda-\lambda_{\rm old}|/\max\{\lambda,\epsilon_{\rm mach}\}$
        \State $\by\gets\by_{\rm new}$, $\lambda_{\rm old}\gets\lambda$
    \Until{$r_y<\epsilon_y$ and $r_\lambda<\epsilon_\lambda$, or $K_{\rm alt}$ alternations have been used}
    \State Recompute $\lambda\gets\texttt{fminbnd}(\operatorname{GCV}(\by,\cdot),10^{-8},10^3)$
    \State Store $(\by^*_\alpha,\lambda^*_\alpha)\gets(\by,\lambda)$
    \State \Return $\operatorname{GCV}(\by^*_\alpha,\lambda^*_\alpha)$ as defined in~\eqref{eq:profiled-gcv-beta}
\EndFunction
\State $\alpha^*\gets\texttt{fminbnd}(\textsc{NestedGCVScore},0,1)$
\State Retrieve $(\by^*,\lambda^*)\gets(\by^*_{\alpha^*},\lambda^*_{\alpha^*})$
\State Compute $\bx^*$ from~\eqref{eq:Tikhonov-solution} using $(\by^*,\lambda^*)$
\State \Return $(\bx^*,\by^*,\lambda^*,\eta_{\rm MAD},\alpha^*)$
\end{algorithmic}
\end{algorithm}

\medskip\noindent\emph{Numerical settings.}
The experiment in Section~\ref{ssec:practical-parameter-choice} uses periodic first-difference image regularization with squared Fourier symbol $4\sin^2(\pi u/256)+4\sin^2(\pi v/256)$, $0\leq u,v<256$, rather than the second-difference regularization used in the earlier experiments.
In Algorithm~\ref{alg:profiled-gcv-beta}, $\epsilon_{\rm mach}$ denotes machine precision.
The $\alpha$ search tolerance is $2\times10^{-3}$, with budgets of 30 evaluations and 30 iterations. The conditional $\lambda$ search tolerance is $10^{-8}$. Alternation allows 200 updates and terminates at relative changes below $10^{-4}$ in $\by$ and $10^{-3}$ in $\lambda$. Each local blur solve allows 50 iterations, with gradient tolerance $10^{-3}/\eta_{\rm MAD}$, initial damping increment $10^{-6}/\eta_{\rm MAD}$, and step tolerance $10^{-12}$. The damping cutoff is $10^8/\eta_{\rm MAD}$; if factorization fails, the fallback direction is $-\eta_{\rm MAD}\nabla \bff_{\beta,\gamma}$. These scaled solver settings preserve the original iteration under deterministic objective scaling.

\section{Additional reconstruction metrics}
\label{app:reconstruction-metrics}
Tables~\ref{tab:additional-gaussian}--\ref{tab:additional-motion} report the remaining metrics for the same 27 runs as Tables~\ref{tab:reconstruction-gaussian}--\ref{tab:reconstruction-motion}. The parameter error and relative residual are defined in Section~\ref{ssec:reconstruction-experiments}; PSNR uses peak intensity one and the unclipped reconstruction. Smaller parameter error and residual and larger PSNR are preferable for their respective criteria, but a small residual alone does not imply accurate recovery. In the tables, PE denotes parameter error, and RR denotes relative residual.

\begin{table}[H]
\centering
\begin{minipage}[t]{0.49\linewidth}
\centering
\small
\caption{\strut Gaussian blur. Runs match Table~\ref{tab:reconstruction-gaussian}.}
\label{tab:additional-gaussian}
\setlength{\tabcolsep}{2pt}
\begin{tabular}{@{}lcrrr@{}}
\toprule
Method & \shortstack{\small Initial\\Guess} & PE & RR & PSNR \\
\midrule
No log-det & 1 & $5.650$ & $3.63\times10^{-4}$ & $20.97$ \\
No log-det & 2 & $5.530$ & $3.63\times10^{-4}$ & $20.97$ \\
No log-det & 3 & $4.960$ & $7.75\times10^{-3}$ & $21.32$ \\
\addlinespace[3pt]
Marg & 1 & $\mathbf{0.123}$ & $9.65\times10^{-3}$ & $\mathbf{23.26}$ \\
Marg & 2 & $\mathbf{0.123}$ & $9.65\times10^{-3}$ & $\mathbf{23.26}$ \\
Marg & 3 & $\mathbf{0.123}$ & $9.65\times10^{-3}$ & $\mathbf{23.26}$ \\
\addlinespace[3pt]
Proposed & 1 & $\mathbf{0.123}$ & $9.70\times10^{-3}$ & $\mathbf{23.26}$ \\
Proposed & 2 & $\mathbf{0.123}$ & $9.70\times10^{-3}$ & $\mathbf{23.26}$ \\
Proposed & 3 & $\mathbf{0.123}$ & $9.70\times10^{-3}$ & $\mathbf{23.26}$ \\
\bottomrule
\end{tabular}
\end{minipage}\hfill
\begin{minipage}[t]{0.49\linewidth}
\centering
\small
\caption{\strut Defocus blur. Runs match Table~\ref{tab:reconstruction-defocus}.}
\label{tab:additional-defocus}
\setlength{\tabcolsep}{2pt}
\begin{tabular}{@{}lcrrr@{}}
\toprule
Method & \shortstack{\small Initial\\Guess} & PE & RR & PSNR \\
\midrule
No log-det & 1 & $2.500$ & $3.76\times10^{-3}$ & $23.33$ \\
No log-det & 2 & $2.500$ & $3.76\times10^{-3}$ & $23.33$ \\
No log-det & 3 & $3.557$ & $2.09\times10^{-2}$ & $5.91$ \\
\addlinespace[3pt]
Marg & 1 & $0.242$ & $9.38\times10^{-3}$ & $25.26$ \\
Marg & 2 & $0.242$ & $9.38\times10^{-3}$ & $25.26$ \\
Marg & 3 & $3.560$ & $2.09\times10^{-2}$ & $5.90$ \\
\addlinespace[3pt]
Proposed & 1 & $\mathbf{0.165}$ & $9.76\times10^{-3}$ & $\mathbf{25.36}$ \\
Proposed & 2 & $\mathbf{0.165}$ & $9.76\times10^{-3}$ & $\mathbf{25.36}$ \\
Proposed & 3 & $3.565$ & $2.09\times10^{-2}$ & $5.89$ \\
\bottomrule
\end{tabular}
\end{minipage}
\end{table}

\begin{table}[H]
\centering
\begin{minipage}[t]{0.55\linewidth}
\centering
\small
\caption{\strut Motion blur. Runs match Table~\ref{tab:reconstruction-motion}.}
\label{tab:additional-motion}
\setlength{\tabcolsep}{2pt}
\begin{tabular}{@{}lcrrr@{}}
\toprule
Method & \shortstack{\small Initial\\Guess} & PE & RR & PSNR \\
\midrule
No log-det & 1 & $10.061$ & $9.41\times10^{-4}$ & $21.64$ \\
No log-det & 2 & $10.447$ & $9.56\times10^{-4}$ & $21.64$ \\
No log-det & 3 & $10.463$ & $9.61\times10^{-4}$ & $21.64$ \\
\addlinespace[3pt]
Marg & 1 & $5.514$ & $6.19\times10^{-3}$ & $20.42$ \\
Marg & 2 & $0.964$ & $8.43\times10^{-3}$ & $23.12$ \\
Marg & 3 & $3.159$ & $1.17\times10^{-2}$ & $19.24$ \\
\addlinespace[3pt]
Proposed & 1 & $5.886$ & $7.37\times10^{-3}$ & $21.61$ \\
Proposed & 2 & $\mathbf{0.055}$ & $9.52\times10^{-3}$ & $\mathbf{24.79}$ \\
Proposed & 3 & $2.010$ & $1.18\times10^{-2}$ & $23.84$ \\
\bottomrule
\end{tabular}
\end{minipage}
\end{table}

\end{document}